\documentclass{LMCS}

\def\dOi{11(1:5)2015}
\lmcsheading%
{\dOi}
{1--28}
{}
{}
{Feb.~21, 2014}
{Jun.~10, 2015}
{}

\ACMCCS{[{\bf Theory of computation}]: Models of
  computation---Probabilistic computation\,/ Quantum computation theory;
  Semantics and reasoning---Program semantics---Categorical semantics}

\newif\ifignore 
\ignorefalse

\newcommand{\auxproof}[1]{
\ifignore\mbox{}\newline
\textbf{PROOF:} \dotfill\newline
{\it #1}\mbox{}\newline
\textbf{ENDPROOF}\dotfill
\fi}

\usepackage[english]{babel}
\usepackage{amsmath}
\usepackage{amssymb}
\usepackage{xspace}
\usepackage{stmaryrd}
\usepackage{url,hyperref}
\usepackage{fancybox}
\usepackage{xypic}
\usepackage{xy}
\xyoption{v2}
\xyoption{curve}
\xyoption{2cell}
\UseTwocells
\usepackage{url}
\usepackage{amssymb}
\usepackage{amstext}
\usepackage{xspace}

\newcommand{\after}{\mathrel{\circ}}
\newcommand{\set}[2]{\{#1\;|\;#2\}}
\newcommand{\setin}[3]{\{#1\in#2\;|\;#3\}}

\newcommand{\all}[2]{\forall{#1}.\,#2}

\newcommand{\lam}[2]{\lambda#1.\,#2}
\newcommand{\lamin}[3]{\lambda#1\in#2.\,#3}
\newcommand{\tuple}[1]{\langle#1\rangle}

\newcommand{\NNO}{\mathbb{N}}
\newcommand{\C}{\mathbb{C}}

\newcommand{\R}{\mathbb{R}}

\newcommand{\U}{\mathcal{U}}
\newcommand{\E}{\mathcal{E}}
\renewcommand{\H}{\mathcal{H}}
\newcommand{\intd}{{\kern.2em}\mathrm{d}{\kern.03em}}
\newcommand{\ket}[1]{\ensuremath{|{\kern.1em}#1{\kern.1em}\rangle}}
\newcommand{\bra}[1]{\langle\,#1\,|}

\newcommand{\unifnorm}[1]{\|#1\|_{\infty}}

\newcommand{\Dst}{\mathcal{D}}

\newcommand{\Giry}{\mathcal{G}}
\newcommand{\Rdn}{\mathcal{R}}
\newcommand{\Ult}{\mathcal{U}}
\newcommand{\Sw}{\ensuremath{\mathcal{\acute{S}}}}
\newcommand{\starfun}[1]{{\ensuremath{\mathcal{C}_{#1}}}}

\newcommand{\Pow}{\mathcal{P}}

\newcommand{\idmap}[1][]{\ensuremath{\mathrm{id}_{#1}}}
\newcommand{\Idmap}[1][]{\ensuremath{\mathrm{Id}_{#1}}}
\newcommand{\op}[1]{#1^{\textrm{op}}}
\newcommand{\SA}{^{\mathrm{sa}}}
\newcommand{\Obj}{\ensuremath{\mathrm{Obj}}}
\newcommand{\orthogonal}{\mathrel{\bot}}
\newcommand{\supp}{\textsl{supp}}

\newcommand{\ev}{\textsl{ev}}
\newcommand{\Hom}{\ensuremath{\mathrm{Hom}}}
\newcommand{\HomMIU}{\ensuremath{\Hom_{\mathrm{MIU}}}}
\newcommand{\HomPU}{\ensuremath{\Hom_{\mathrm{PU}}}}
\newcommand{\Cont}{\textsl{Cont}}

\newcommand{\Stat}{\ensuremath{\mathrm{Stat}}}
\newcommand{\MStat}{\ensuremath{\mathrm{MStat}}}

\newcommand{\cat}[1]{\ensuremath{\textbf{#1}}\xspace}

\newcommand{\Sets}{\ensuremath{\textbf{Sets}}\xspace}
\newcommand{\FinSets}{\ensuremath{\textbf{FinSets}}\xspace}

\newcommand{\EA}{\ensuremath{\textbf{EA}}\xspace}

\newcommand{\CH}{\ensuremath{\textbf{CH}}\xspace}

\newcommand{\CCL}{\ensuremath{\cat{CCLcvx}}}

\newcommand{\EMod}{\ensuremath{\textbf{EMod}}\xspace}

\newcommand{\Conv}{\ensuremath{\textbf{Conv}}\xspace}

\newcommand{\Vect}{\ensuremath{\textbf{Vect}}\xspace}

\newcommand{\CstarMap}[1]{\ensuremath{\textbf{Cstar}_{#1}}\xspace}
\newcommand{\CstarPU}{\CstarMap{\textrm{PU}}}
\newcommand{\CstarCPU}{\CstarMap{\textrm{cPU}}}
\newcommand{\CstarMIU}{\CstarMap{\textrm{MIU}}}

\newcommand{\CCstarMap}[1]{\ensuremath{\textbf{CCstar}_{#1}}\xspace}
\newcommand{\CCstarPU}{\CCstarMap{\textrm{PU}}}
\newcommand{\CCstarMIU}{\CCstarMap{\textrm{MIU}}}

\newcommand{\FdCstarMap}[1]{\ensuremath{\textbf{FdCstar}_{#1}}\xspace}

\newcommand{\FdCCstarMap}[1]{\ensuremath{\textbf{FdCCstar}_{#1}}\xspace}
\newcommand{\FdCCstarPU}{\FdCCstarMap{\textrm{PU}}}
\newcommand{\FdCCstarMIU}{\FdCCstarMap{\textrm{MIU}}}

\newcommand{\scalar}{\mathrel{\bullet}}

\newcommand{\inv}{\mathop{\rlap{\raisebox{.3ex}{${\kern.6ex}\cdot$}}-}}

\newcommand{\id}{\ensuremath{\mathrm{id}}}

\newcommand{\Kl}{\mathcal{K}{\kern-.2ex}\ell}

\newcommand{\KlN}{\Kl_{\NNO}}
\newcommand{\EM}{\mathcal{E}{\kern-.2ex}\mathcal{M}}

\newcommand{\Ef}{\ensuremath{\mathcal{E}{\kern-.5ex}f}}
\newcommand{\DM}{\ensuremath{\mathcal{D}{\kern-.85ex}\mathcal{M}}}
\newcommand{\klafter}{\mathrel{\raisebox{.15em}{$\scriptscriptstyle\odot$}}}
\newcommand{\Aff}{\ensuremath{\mathcal{A}}}

\newcommand{\conglongrightarrow}{\mathrel{\smash{\stackrel{
           \raisebox{.5ex}{$\scriptstyle\cong$}}{
           \raisebox{0ex}[0ex][0ex]{$\longrightarrow$}}}}}

\makeatother
 \newdir{ >}{{}*!/-7.5pt/@{>}}
 \newdir{|>}{!/4.5pt/@{|}*:(1,-.2)@^{>}*:(1,+.2)@_{>}}
 \newdir{ |>}{{}*!/-3pt/@{|}*!/-7.5pt/:(1,-.2)@^{>}*!/-7.5pt/:(1,+.2)@_{>}}
\newcommand{\xyline}[2][]{\ensuremath{\smash{\xymatrix@1#1{#2}}}}
\newcommand{\xyinline}[2][]{\ensuremath{\smash{\xymatrix@1#1{#2}}}^{\rule[8.5pt]{0pt}{0pt}}}
\makeatletter

\begin{document}

\title[Probabilistic Gelfand Duality]{From Kleisli categories to commutative $C^*$-algebras: \\
   Probabilistic Gelfand Duality}

\author[R.~Furber]{Robert Furber}
\address{Institute for Computing and Information Sciences (iCIS) \\
Radboud University Nijmegen, The Netherlands.}
\email{\{r.furber,bart\}@cs.ru.nl}

\author{Bart Jacobs}
\address{\vspace{-18 pt}}



\begin{abstract}
$C^*$-algebras form rather general and rich mathematical structures
  that can be studied with different morphisms (preserving
  multiplication, or not), and with different properties (commutative,
  or not). These various options can be used to incorporate various
  styles of computation (set-theoretic, probabilistic, quantum) inside
  categories of $C^*$-algebras. At first, this paper concentrates on
  the commutative case and shows that there are functors from several
  Kleisli categories, of monads that are relevant to model
  probabilistic computations, to categories of $C^*$-algebras. This
  yields a new probabilistic version of Gelfand duality, involving the
  ``Radon'' monad on the category of compact Hausdorff spaces. We then
  show that the state space functor from $C^*$-algebras to
  Eilenberg-Moore algebras of the Radon monad is full and
  faithful. This allows us to obtain an appropriately commuting
  state-and-effect triangle for $C^*$-algebras.
\end{abstract}

\maketitle

\section{Introduction}
There are several notions of computation. We have the classical notion
of computation, probabilistic computation, where a computer may make
random choices, and quantum computation, which uses quantum mechanical
interference and measurement. Normally we would consider classical
computation to be done on sets, probabilistic computation on spaces
with a measure, and quantum computation on Hilbert spaces. We can
instead use categories with $C^*$-algebras as objects and a choice of
either *-homomorphisms (called MIU-map below) or positive unital maps
as the morphisms. We note at this point that positive unital maps coincide with \emph{completely} positive unital maps if either the domain or codomain of a map is a commutative $C^*$-algebra, but not in general. The general outline is represented in this table.
\newpage
\begin{center}\renewcommand{\tabcolsep}{.7em}
\begin{tabular}{c||c|c|c}
 & \textbf{set-theoretic} & 
   \textbf{probabilistic} & \textbf{quantum}
\vrule height5mm depth3mm width0mm \\
\hline\hline
$C^*$-algebras & commutative & commutative & non-commutative 
\vrule height5mm depth3mm width0mm \\
\hline
maps preserve & 
\begin{tabular}{c} multiplication \\[-.6em] involution \\[-.6em] 
   unit \end{tabular} &
\begin{tabular}{c} positivity \\[-.6em] unit \end{tabular} &
\begin{tabular}{c} positivity \\[-.6em] unit \end{tabular} \\
\hline
maps abbreviation & MIU & PU & PU 
\vrule height5mm depth3mm width0mm \\
\end{tabular}
\end{center}

\noindent While the quantum case is an important source of motivation,
we will deal more with the classical and probabilistic cases in
this article. In particular, we will relate the alternative method of
representing probabilistic computation, using monads, to the
$C^*$-algebraic approach.

In recent years the methods and tools of category theory have been
applied to Hilbert spaces --- see \textit{e.g.}~\cite{AbramskyC09} and
the references there --- and also to $C^{*}$-algebras, see for
instance~\cite{Pelletier93,DualityMarkov}.  In this paper we show that
clearly distinguishing different types of homomorphisms of
$C^*$-algebras already brings quite some clarity.  Moreover, we
demonstrate the relevance of monads (and their Kleisli and
Eilenberg-Moore categories) in this field. The aforementioned
paper~\cite{Pelletier93} concerns itself with only the *-homomorphisms
(\textit{i.e.}~with the MIU-maps in our terminology). 

The main results of the paper can be summarised as follows. The
well-known finite (`baby') version of Gelfand duality involves an
equivalence between on the one hand the category of finite sets (and
all functions between them), and on the other hand the opposite of the
category of finite-dimensional commutative $C^*$-algebras with
MIU-maps (*-homomorphisms) between them. Diagrammatically:
$$\xymatrix{
\FinSets\ar[r]^-{\simeq} & \op{\big(\FdCCstarMIU\big)}
}$$

\noindent Our first observation is that if we generalise from MIU to
PU (positive unital) maps we get an equivalence:
$$\xymatrix{
\KlN(\Dst)\ar[r]^-{\simeq} & \op{\big(\FdCCstarPU\big)}
}$$

\noindent where $\Dst$ is the distribution monad on $\Sets$, and
$\KlN(\Dst)$ is the Kleisli category of this monad, but with objects
restricted to natural numbers. This shows that the category
$\FdCCstarPU$ is the Lawvere theory of the distribution monad. Details
are in Section~\ref{DiscProbSec}.

The main contribution of the paper lies in a generalisation of the
latter equivalence beyond the finite case, which can be summarised
in a diagram:
\begin{equation}
\label{OverviewDiag}
\vcenter{\xymatrix{
\CH\ar[rr]^-{\simeq}_-{\text{Gelfand}}\ar@/_2.5ex/[d]\ar@(u,l)_{\Rdn} & & 
   \op{\big(\CCstarMIU\big)}\ar@{^(->}[d] \\
\Kl(\Rdn)\ar[rr]^-{\simeq}_-{\text{new}}\ar[u]^{\dashv} & & 
   \op{\big(\CCstarPU\big)}
}}
\end{equation}

\noindent At the top of this diagram we have the classical Gelfand
duality between the category $\CH$ of compact Hausdorff spaces and the
(opposite of the) category of commutative $C^*$-algebras with
MIU-maps. Again, the generalisation to the computationally more
interesting PU-maps involves a duality with a Kleisli category, namely
the Kleisli category $\Kl(\Rdn)$ of what we call the Radon monad
$\Rdn$ on compact Hausdorff spaces. Elements of $\Rdn(X)$ can be
described as so-called Radon probability measures, also known as inner
regular probability measures (see~\cite{Rudin87}).

In the end, in Diagram~\eqref{KlRSETDiag} we show how the Kleisli
category of the Radon monad gives rise to a `state-and-effect'
triangle that combines Kleisli computations for the Radon monad and
their associated predicate transformers and state transformers. These
predicate and state transformers correspond to the Heisenberg and
Schr\"odinger picture, respectively.

Incidentally, the adjunction on the left in
Diagram~\eqref{OverviewDiag} can be transferred to the right, and then
yields a right adjoint to the inclusion $\CCstarMIU \hookrightarrow
\CCstarPU$. In~\cite{BramWesterbaan14} it is shown that such a right
adjoint also exists in the general non-commutative case.

Giry~\cite[I.4]{Giry82} described how we can consider a stochastic
process as being a diagram in the Kleisli category of the Giry monad
on measure spaces. By using the Radon monad $\Rdn$ on compact spaces
instead, we can get a different category of stochastic processes on
compact spaces as diagrams in the (opposite of the) category of
\emph{commutative} $C^*$-algebras with PU-maps. This allows the
quantum generalization to taking diagrams in the category of
\emph{non-commutative} $C^*$-algebras, or by considering diagrams in
the category $\EM(\Rdn)$ of Eilenberg-Moore algebras of the Radon
monad $\Rdn$, in which the category of $C^*$-algebras faithfully
embeds. The relationship to quantum computation is that $B(\H)$, the
algebra of all bounded operators on a Hilbert space $\H$, is a
$C^*$-algebra, and for every $C^*$-algebra $A$, there is a Hilbert
space $\H$ such that $A$ is isomorphic to a norm-closed *-subalgebra
of $B(\H)$. Unitary maps $U \colon \H \rightarrow \H$ define MIU maps
$a \mapsto U^*aU \colon B(\H) \rightarrow B(\H)$. The category of
$C^*$-algebras allows us to represent measurement with maps from a
commutative $C^*$-algebra to $B(\H)$. We can also represent composite
systems that are partly quantum and partly classical. Girard also used
certain special $C^*$-algebras, von Neumann algebras, for his Geometry
of Interaction~\cite{Girard2011}.

\section{Preliminaries on \texorpdfstring{$C^*$}{C*}-algebras}\label{PrelimSec}

We write $\Vect = \Vect_{\C}$ for the category of vector spaces over
the complex numbers $\C$. This category has direct product $V\oplus
W$, forming a biproduct (both a product and a coproduct) and tensors
$V\otimes W$, which distribute over $\oplus$. The tensor unit is the
space $\C$ of complex numbers. The unit for $\oplus$ is the singleton
(null) space $0$. We write $\overline{V}$ for the vector space with
the same vectors/elements as $V$, but with conjugate scalar product:
$z\scalar_{\overline{V}} v = \overline{z} \scalar_{V} v$. This makes
$\Vect$ an involutive category, see~\cite{Jacobs12f}.

A \emph{*-algebra} is an involutive monoid $A$ in the category
$\Vect$.  Thus, $A$ is itself a vector space, carries a multiplication
$\cdot \colon A\otimes A \rightarrow A$, linear in each argument, and
has a unit $1\in A$. Moreover, there is an involution map
$(-)^{*} \colon \overline{A} \rightarrow A$, preserving $0$ and $+$
and satisfying:
$$\begin{array}{rclcrclcrclcrcl}
1^{*} & = & 1
& \quad &
(x\cdot y)^{*} & = & y^{*}\cdot x^{*}
& \quad &
x^{**} & = & x
& \quad &
(z \scalar x)^{*}
& = &
\overline{z} \scalar x^{*}.
\end{array}$$

\noindent Here we have written a fat dot $\scalar$ for scalar
multiplication, to distinguish it from the algebra's multiplication
$\cdot$. For $z=a+bi\in\C$ we have the conjugate $\overline{z} =
a-bi$. Often we omit the multiplication dot $\cdot$ and simply write
$xy$ for $x\cdot y$. Similarly, the scalar multiplication $\scalar$ is
often omitted. We then rely on the context to distinguish the two
multiplications.

A \emph{$C^*$-algebra} is a *-algebra $A$ with a norm $\|-\|\colon A
\rightarrow \R_{\geq 0}$ in which it is complete, satisfying the conditions
$\|x\| = 0$ iff $x=0$ and:
$$\begin{array}{rclcrcl}
\|x+y\| & \leq & \|x\| + \|y\|
& \qquad &
\|z\scalar x\| & = & |z| \cdot \|x\| \\
\|x\cdot y\| & \leq & \|x\| \cdot \|y\|
& & 
\|x^{*}\cdot x\| & = & \|x\|^{2}.
\end{array}$$

\noindent The last equation $\|x^{*}\cdot x\| = \|x\|^{2}$, is the
\emph{$C^*$-identity} and distinguishes $C^*$-algebras from Banach
*-algebras. We remark at this point that a Banach *-algebra admits at
most one norm satisfying the $C^*$-identity. The reason for this is
that the spectral radius $r(x)$ is definable in terms of the ring
structure of the algebra, and for self-adjoint elements $r(x) = \|x\|$
\cite[Proposition 4.1.1 (a)]{KadisonR83}. If $x$ is an arbitrary
element, $x^* \cdot x$ is self-adjoint, so $r(x^* \cdot x) = \|x^*
\cdot x\| = \|x\|^2$.  In the current setting, each $C^*$-algebra is
unital, \textit{i.e.}~has a (multiplicative) unit $1$. A consequence
of the axioms above is that $\|1\| = 1$ unless the $C^*$-algebra is
the unique one in which $0 = 1$. A $C^*$-algebra is called
\emph{commutative} if its multiplication is commutative, and
\emph{finite-dimensional} is it has finite dimension when considered
as a vector space.

An element $x$ in a $C^*$-algebra $A$ is called \emph{positive} if it
can be written in the form $x = y^{*}\cdot y$.  We write
$A^{+}\subseteq A$ for the subset of positive elements in $A$.  This
subset is a cone, which is to say it is closed under addition and
scalar multiplication with positive real numbers. The multiplication
$x\cdot y$ of two positive elements need not be positive in general
(think of matrices). The square $x^{2} = x\cdot x$ of a self-adjoint
element $x = x^{*}$, however, is obviously positive. In a
\emph{commutative} $C^*$-algebra the positive elements are closed
under multiplication. A cone $A^+$ in a vector space defines a partial
order as follows.
\begin{equation}
\label{ConePosetDef}
\begin{array}{rcl}
x \leq y 
& \Longleftrightarrow &
y - x \in A^{+}.
\end{array}
\end{equation}
\noindent This is defines an order on every $C^*$-algebra.

\auxproof{
Positive elements are self-adjoint: $x^{*} = (y^{*}\cdot y)^{*} =
y^{*}\cdot y^{**} = y^{*}\cdot y = x$.

If $x\in A$ is positive, say $x = u^{*}u$, then scalar multiplication
$rx\in A$ with a positive $r\in\R_{\geq 0}$ yields a positive element,
since:
$$\begin{array}{rcl}
(\sqrt{r}\scalar u)^{*}\cdot (\sqrt{r}\scalar u)
& = &
(\overline{\sqrt{r}} \scalar u^{*}) \cdot (\sqrt{r}\scalar u) \\
& = &
(\sqrt{r} \scalar u^{*}) \cdot (\sqrt{r}\scalar u) \\
& = &
\sqrt{r} \scalar ((\sqrt{r}\scalar u^{*})\cdot  u) \\
& = &
\sqrt{r} \scalar (\sqrt{r}\scalar (u^{*}\cdot  u)) \\
& = &
(\sqrt{r} \cdot \sqrt{r})\scalar x \\
& = &
r\scalar x.
\end{array}$$

If $x = u^{*}u$ then $x^{2}$ is also positive, since:
$$x^{2} = u^{*}uu^{*}u = u^{*}u^{**}u^{*}u = (u^{*}u)^{*}u^{*}u.$$

Suppose $x = u^{*}u$ and $y = v^{*}v$, then, in the commutative case:
$$(xy)^{*} = y^{*}x^{*} = yx = xy = u^{*}uv^{*}v = u^{*}v^{*}uv
   = (vu)^{*}vu.$$

Some more useful facts.
\begin{enumerate}
\item $\|x^{*}\| = \|x\|$ and $\|1\| = 1$. A proof of the first
  equation can be found in~\cite[Lemma 1.1.6]{Sakai71}: $\|x^{2}\|^{2}
  = \|x^{*}x\| \leq \|x^{*}\| \|x\|$, so $\|x\| \leq \|x^{*}\|$, for
  each $x$.  But then also $\|x^{*}\| \leq \|x^{**}\| = \|x\|$. Hence
  $\|x^{*}\| = \|x\|$. Further, $\|1\| = \|1^{*}1\| = \|1\|^{2}$, and
  since $1\neq 0$ we get $\|1\|\neq 0$ and thus $\|1\| = 1$.

\item $x \leq \|x\| \scalar 1$, for positive $x$, see Warner, ``NB''
before lemma 1.29.

\item For positive elements $x,y$, if $x\leq y$, then $\|x\| \leq \|y\|$,
see Warner, lemma 1.25

\item positive elements $x$ have unique $n$-th square roots
  $x^{\frac{1}{n}}$, see~\cite[Prop. 1.4.1]{Sakai71}. Hence for a
  self-adjoint $x$, put $|x| = \sqrt{x^{2}}$.

\item self-adjoint elements carry a partial order, which is preserved
by positive maps.

\end{enumerate}

}

There are mainly two options when it comes to maps between
$C^*$-algebras. The difference between them plays an important role in
this paper.

\begin{defi}
\label{CstarCatDef}
We define two categories $\CstarMIU$ and $\CstarPU$ with $C^*$-algebras
as objects, but with different morphisms.
\begin{enumerate}
\item A morphism $f\colon A \rightarrow B$ in $\CstarMIU$ is a linear
  map preserving multiplication (M), involution (I), and unit
  (U). Explicitly, this means for all $x,y\in A$,
$$\begin{array}{rclcrclcrcl}
f(x\cdot y) & = & f(x)\cdot f(y)
& \qquad & 
f(x^{*}) & = & f(x)^{*}
& \qquad & 
f(1) & = & 1.
\end{array}$$

\noindent Often such ``MIU'' maps are called *-homomorphisms.

\item A morphism $f\colon A \rightarrow B$ in $\CstarPU$ is a linear
  map that preserves positive elements and the unit. This means that
  $f$ restricts to a function $A^{+} \rightarrow
  B^{+}$. Alternatively, for each $x\in A$ there is an $y\in B$ with
  $f(x^{*}x) = y^{*}y$. \qed
\end{enumerate}

\noindent For both $X = \textrm{MIU}$ and $X=\textrm{PU}$ there are
obvious full subcategories of commutative and/or finite-dimensional
$C^*$-algebras, as described in:
$$\xymatrix@R-1.3pc@C+1pc{
& \CCstarMap{X}\ar@{^(->}[dr] & \\
\FdCCstarMap{X}\ar@{^(->}[ur]\ar@{^(->}[dr] & & \CstarMap{X} \\
& \FdCstarMap{X}\ar@{^(->}[ur]
}$$
\end{defi}

Clearly, each ``MIU'' map is also a ``PU'' map, so that we have
inclusions $\CstarMIU \hookrightarrow \CstarPU$, also for the various
subcategories. A map that preserves positive elements is called
positive itself; and a unit preserving map is called unital.  Positive
unital maps are the natural notion of morphism between order unit
spaces and Riesz spaces.

For a category $\cat{B}$ one often writes $\cat{B}(X,Y)$ or
$\Hom(X,Y)$ for the ``homset'' of morphisms $X\rightarrow Y$ in
$\cat{B}$. For $C^*$-algebras $A,B$ we write $\HomMIU(A,B) =
\CstarMIU(A,B)$ and $\HomPU(A,B) = \CstarPU(A,B)$ for the homsets of
MIU- and PU-maps. For the special case where $B$ is the algebra $\C$
of complex numbers we define sets of ``states'' and of
``multiplicative states'' as:
$$\begin{array}{rclcrcl}
\Stat(A)
& = &
\HomPU(A,\C)
& \qquad\mbox{and}\qquad &
\MStat(A)
& = &
\HomMIU(A,\C).
\end{array}$$

There is also the commonly used notion of completely positive maps,
which is a stronger condition than positivity but weaker than being
MIU. These maps are important when defining the tensor of
$C^*$-algebras as a functor, as the tensor of positive maps need not
be positive. They are also widely considered to represent the
physically realizable transformations. Positive, but non-completely
positive maps of $C^*$-algebras also have their uses, as entanglement
witnesses for example~\cite[theorem 2]{horodecki}. Since we mainly
consider the commutative case, where positive and completely positive
coincide, we do not consider the category of $C^*$-algebras with
completely positive maps any further in this paper. However, since a
completely positive unital map is what is known as a channel in
quantum information, then theorem \ref{KlRToStarThm} shows that every
channel in Mislove's sense \cite{Mislove2012} is a channel in this
sense.

We collect some basic (standard) properties of PU-morphisms between
$C^*$-algebras (see \textit{e.g.}~\cite{Sakai71,Arveson81}).

\begin{lem}
\label{CstarMapLem}
A PU-map, \textit{i.e.}~a morphism in the category $\CstarPU$,
commutes with involution $(-)^{*}$, and preserves the partial order
$\leq$ given by \eqref{ConePosetDef}.

Moreover, a PU-map $f$ satisfies $\|f(x)\| \leq 4\|x\|$, so that 
$\|f(x) - f(y)\| \leq 4\|x-y\|$, making $f$ continuous.
\end{lem}

\begin{proof}
An element $x$ is called self-adjoint if $x^{*} = x$. Each
self-adjoint $x$ can be written uniquely as a difference $x = x_{p} - x_{n}$ of
positive elements $x_{p}, x_{n}$, with $x_{p}x_{n} = x_{n}x_{p} = 0$ and $\|x_{p}\|, \|x_{n}\| \leq \|x\|$, see~\cite[Proposition 4.2.3 (iii)]{KadisonR83}; as a result $f(x^{*}) = f(x) = f(x)^{*}$, for a PU-map $f$. Next, an arbitrary element $y$ can be
written uniquely as $y = y_{r} + iy_{i}$ for self-adjoint elements $y_{r} =
\frac{1}{2}(y + y^{*}), y_{i} = \frac{1}{2i}(y-y^{*})$, so that
$\|y_{r}\|, \|y_{i}\| \leq \|y\|$. Then $f(y^{*}) =
f(y)^{*}$. Preservation of the order is trivial.

For positive $x$ we have $x \leq \|x\| \scalar 1$, and thus $f(x) \leq
\|x\|\scalar 1$, which gives $\|f(x)\| \leq \|x\|$. An arbitrary
element $x$ can be written as linear combination of four positive
elements $x_{i}$, as in $x = x_{1} - x_{2} + ix_{3} - ix_{4}$, with
$\|x_{i}\| \leq \|x\|$. Finally, $\|f(x)\| = \|f(x_{1}) - f(x_{2}) +
if(x_{3}) - if(x_{4})\| \leq \sum_{i}\|f(x_{i})\| \leq
\sum_{i}\|x_{i}\| \leq 4\|x\|$. 
\end{proof}

\auxproof{
That $x$ can be written as difference of positive elements may
be found in~\cite[Defn. 1.4.3]{Sakai71}. Explicitly, in terms of
the absolute value $|x| = \sqrt{x^{2}}$ one can define $x_{p} = 
\frac{1}{2}(|x| + x)$ and $x_{n} = \frac{1}{2}(|x|-x)$. Then $x_{p}
\cdot x_{n} = 0$. And $x_{p}, x_{n}$ are unique with these properties,
see Warner Lemma 1.28.

Since $f$ preserves positive 
elements, we have:
$$\begin{array}{rcccccl}
f(x)^{*}
& = &
f(x_{p})^{*} - f(x_{n})^{*}
& = &
f(x_{p}) - f(x_{n}) 
& = &
f(x).
\end{array}$$

Take $y_{r} = \frac{1}{2}(y + y^{*})$ and $y_{i} = \frac{1}{2}(iy^{*} - iy)$,
then $y_{r}, y_{i}$ are self-adjoint, and $y = y_{r} + iy_{i}$:
$$\begin{array}{rcl}
y_{r}^{*}
& = &
\frac{1}{2}(y^{*} + y^{**}) \\
& = &
\frac{1}{2}(y + y^{*}) \\
& = &
y_{r} \\
y_{i}^{*}
& = &
\frac{1}{2}(-iy^{**} + iy^{*}) \\
& = &
\frac{1}{2}(iy^{*} - iy) \\
& = &
y_{i} \\
y_{r} + iy_{i}
& = &
\frac{1}{2}(y + y^{*}) + \frac{1}{2}i(iy^{*} - iy) \\
& = &
\frac{1}{2}y + \frac{1}{2}y^{*} - \frac{1}{2}y^{*} + \frac{1}{2}y \\
& = &
y.
\end{array}$$

\noindent And:
$$\begin{array}{rcl}
f(y^{*})
& = &
f(y_{r}^{*} - iy_{i}^{*}) \\
& = &
f(y_{r} - iy_{i}) \\
& = &
f(y_{r}) - if(y_{i}) \\
& = &
f(y_{r})^{*} - if(y_{i})^{*} \\
& = &
\big(f(y_{r}) + f(iy_{i})\big)^{*} \\
& = &
f(y_{r} + iy_{i})^{*} \\
& = &
f(y)^{*}.
\end{array}$$
}

In fact, it can be shown that $\|f(x)\| \leq \|x\|$ for all $x$, not
just positive $x$, reducing the constant 4 in the inequality above to
1 (see~\cite[corollary 1]{RussoDye66}). But this sharpening is not
needed here.

We next recall two famous adjunctions involving compact Hausdorff
spaces.  The first one is due to Manes~\cite{Manes69} and describes
compact Hausdorff spaces as monadic over $\Sets$, via the ultrafilter
monad. The second one is known as Gelfand duality, relating compact
Hausdorff spaces and \emph{commutative} $C^*$-algebras. Notice that
this result involves the ``MIU'' maps.

\begin{thm}
\label{GelfandThm}
Let $\CH$ be the category of compact Hausdorff spaces, with continuous
maps between them. There are two fundamental adjunctions:
$$\xymatrix@R-1.5pc{
\CH\ar@/^2ex/[dd]^-{\text{forget}}
& \qquad &
\CH\ar@/_2ex/[dd]_{C} \\
\dashv & & \simeq \\
\Sets\ar@/^2ex/[uu]^{\Ult}
& & 
\op{(\CCstarMIU)}\ar@/_2ex/[uu]_{\MStat}
}$$

\noindent On the left the functor $\Ult$ sends a set $X$ to the
ultrafilters on the powerset $\Pow(X)$. And on the right the
equivalence of categories is given by sending a compact Hausdorff
space $X$ to the commutative $C^*$-algebra $C(X) = \Cont(X,\C)$ of
continuous functions $X \rightarrow \C$. The ``weak-* topology'' on
states will be discussed below. \qed
\end{thm}

The multiplicative states on a commutative $C^*$-algebra can
equivalently be described as maximal ideals, or also as so-called pure
states (see below).

\begin{cor}
\label{FinGelfandCor}
For each finite-dimensional commutative $C^*$-algebra $A$ there is
an $n\in\NNO$ with $A \cong \C^{n}$ in $\FdCCstarMIU$.
\end{cor}

\begin{proof}
By the previous theorem there is a compact Hausdorff space $X$ such
that $A$ is MIU-isomorphic to the algebra of continuous maps $X
\rightarrow \C$.  This $X$ must be finite, and since a finite
Hausdorff space is discrete, all maps $X\rightarrow \C$ are
continuous. Let $n\in\NNO$ be the number of elements in $X$; then we
have an isomorphism $A \cong \C^{n}$.
\end{proof}

As we can already see in the above theorem, it is the \emph{opposite}
of a category of $C^*$-algebras that provides the most natural setting
for computations. This is in line with what is often called the
Heisenberg picture. In a logical setting it corresponds to
computation of weakest preconditions, going backwards. The situation
may be compared to the category of complete Heyting algebras, which is
most usefully known in opposite form, as the category of locales,
see~\cite{Johnstone82}.

{\sloppy The set of states $\Stat(A) = \HomPU(A,\C)$ now can be equipped with the
weak-* topology, defined as the coarsest (smallest) topology in which
all evaluation maps $\ev_{x} = \lam{s}{s(x)} \colon \HomPU(A,\C)
\rightarrow \C$, for $x\in A$, are continuous. We introduce the category $\CCL$, which first appeared in \cite{Swirszcz74}, in order to extend $\Stat$ to a functor. }

The category $\CCL$ has as its objects compact convex subsets of
(Hausdorff) locally convex vector spaces. More accurately, the objects
are pairs $(V, X)$ where $V$ is a (Hausdorff) locally convex space,
and $X$ is a compact convex subset of $V$. The maps $(V,X) \rightarrow
(W,Y)$ are continuous, affine maps $X \rightarrow Y$. Note that if
$(V,X)$ and $(W,Y)$ are isomorphic, while $X$ is necessarily
homeomorphic to $Y$, $V$ need not bear any particular relation to $W$
at all. We can see $\CCL$ forms a category, as identity maps are
affine and continuous and both of these attributes of a map are
preserved under composition. We remark at this point that we have a
forgetful functor $U \colon \CCL \rightarrow \CH$, taking the
underlying compact Hausdorff space of $X$.

\begin{prop}
\label{StatesProp}
For each $C^*$-algebra $A$, the set of states $\Stat(A) =
\HomPU(A,\C)$ is convex, and is a compact Hausdorff subspace of the
dual space of $A$ given the weak-* topology. Each PU-map $f\colon A
\rightarrow B$ yields an affine continuous function $\Stat(f) =
(-)\after f \colon \Stat(B) \rightarrow \Stat(A)$. This defines a
functor $\Stat \colon \op{(\CstarPU)} \rightarrow \CCL$.

\end{prop}

We recall that a function (between convex sets) is called
\emph{affine} if it preserves convex sums. We will see shortly that
such affine maps are homomorphisms of Eilenberg-Moore algebras for the
distribution monad $\Dst$.

\begin{proof}
For each finite collection $h_{i}\in\HomPU(A,\C)$ with $r_{i}\in[0,1]$
satisfying $\sum_{i}r_{i}=1$, the function $h = \sum_{i}r_{i}h_{i}$ is
again a state.  Moreover, such convex sums are preserved by
precomposition, making the maps $(-) \after f$ affine.
\auxproof{
If states $h,k\in\HomPU(A,\C)$ are not equal, say $h(x)\neq k(x)$,
then either the real or imaginary parts of $h(x)$ and $k(x)$
differ. Let's consider the former. Then there is a real number $r$
with $\mathrm{re}(h(x)) < r < \mathrm{re}(k(x))$. The two open subsets
$\ev_{x}^{-1}(\setin{z}{\C}{\mathrm{re}(z) < r})$ and
$\ev_{x}^{-1}(\setin{z}{\C}{\mathrm{re}(z) > r})$ then separate $h,k$.
Hence $\HomPU(A,\C)$ is Hausdorff.}

The fact that the dual space of $A$, given the weak-* topology, is a
locally convex space is standard, and only uses that $A$ is a Banach
space \cite[Example 1.8]{Conway90}. This implies that the space of states is Hausdorff. The space
of states is closed since because the positive cone in a $C^*$-algebra is closed \cite[Proposition 2.4.5 (i)]{KadisonR83}\cite[Proposition 1.6.1]{Dixmier77} and the set of linear functionals such that $\phi(1) = 1$ is weak-* closed, and the set of states is the intersection of the two. The space of states is also
bounded as each state has norm $1$. Therefore the state space is a
closed and bounded and hence compact by the Banach-Alaoglu Theorem.
\auxproof{ Proof that a weak-* convergent net of states converges to a
  state:

Let $\phi_\alpha$ be a net of states that converges to $\phi$ in the
weak-* topology, i.e. for all $a \in A$, $\phi_\alpha(a)$ converges to
$\phi(a)$ in $\C$.
\begin{itemize}
\item Preservation of unit: We have that $\phi_\alpha(1) = 1$ for all $\alpha$, so $\phi(1) = 1$.
\item Positivity: Let $p \in A$ be a positive element. We have that $\phi_\alpha(p) \in [0,\infty)$ for all $\alpha$, and $[0,\infty)$ is a closed subset of $\C$, and so $\phi(p) \in [0,\infty)$ as well.
\end{itemize}
}

Precomposition $(-)\after f$ is continuous, since for $x\in A$ and
$U\subseteq \C$ open we get an open subset $\big((-)\after
f\big)^{-1}(\ev_{x}^{-1}(U)) = \set{h}{\ev_{x}(h\after f) \in U} =
\ev_{f(x)}^{-1}(U)$. 

Precomposition with the identity map gives the same state again, so
$\Stat$ preserves identity maps. Since composition of PU-maps is
associative, $\Stat$ preserves composition, and hence is a
functor.
\end{proof}


\subsection{Effect modules}\label{EModSubsec}

Effect algebras have been introduced in mathematical
physics~\cite{FoulisB94}, in the investigation of quantum probability,
see~\cite{DvurecenskijP00} for an overview.  An \emph{effect algebra}
is a partial commutative monoid $(M, 0, \ovee)$ with an
orthocomplement $(-)^{\perp}$. One writes $x\orthogonal y$ if $x\ovee
y$ is defined. The formulation of the commutativity and associativity
requirements is a bit involved, but essentially straightforward.  The
orthocomplement satisfies $x^{\perp\perp} = x$ and $x\ovee x^{\perp} =
1$, where $1 = 0^{\perp}$. There is always a partial order, given by
$x\leq y$ iff $x \ovee z = y$, for some $z$.  The main example is the
unit interval $[0,1] \subseteq \mathbb{R}$, where addition $+$ is
obviously partial, commutative, associative, and has $0$ as unit;
moreover, the orthocomplement is $r^{\perp} = 1 -r$. 
We write $\EA$ for the category of effect algebras, with morphism
preserving $\ovee$ and $1$ --- and thus all other structure.

For each set $X$, the set $[0,1]^{X}$ of fuzzy predicates on $X$ is an
effect algebra, via pointwise operations. Each Boolean algebra $B$ is
an effect algebra with $x\orthogonal y$ iff $x\wedge y = \bot$; then
$x\ovee y = x\vee y$. In a quantum setting, the main example is the
set of effects $\Ef(H) = \set{E\colon H \rightarrow H}{0 \leq E \leq
  I}$ on a Hilbert space $H$, see
\textit{e.g.}~\cite{DvurecenskijP00,HeinosaariZ12}.

An \emph{effect module} is an ``effect'' version of a vector space.
It involves an effect algebra $M$ with a scalar multiplication
$s\scalar x\in M$, where $s\in [0,1]$ and $x\in M$. This scalar
multiplication is required to be a suitable homomorphism in each
variable separately. The algebras $[0,1]^{X}$ and $\Ef(H)$ are clearly
such effect modules. Maps in $\EMod$ are $\EA$ maps that are
additionally required to commute with scalar multiplication.

For a $C^*$-algebra $A$ the subset $A^{+} \hookrightarrow A$ of
positive elements carries a partial order $\leq$ defined on
self-adjoint elements in \eqref{ConePosetDef}. We write $[0,1]_{A} \subseteq A^{+} \subseteq
A$ for the subset of positive elements below the unit. The elements in
$[0,1]_{A}$ will be called effects (or sometimes also: predicates).
For instance, for the $C^*$-algebra $B(\H)$ of bounded operators on a
Hilbert space $\H$ the unit interval $[0,1]_{B(\H)} \subseteq B(H)$
contains the effects $\Ef(\H) = \setin{A}{B(\H)}{0 \leq A \leq
  \idmap}$ on $\H$.

We claim that $[0,1]_{A}$ is an \emph{effect algebra} and carries a
$[0,1] \subseteq \R$ scalar multiplication, thus making it an
\emph{effect module}.
\begin{itemize}
\item Since $A$ with $0, +$ is a partially ordered Abelian group,
  $[0,1]_{A}$ is a so-called interval effect algebra, with
  $x\orthogonal y$ iff $x+y\leq 1$, and in that case $x\ovee y =
  x+y$. The orthocomplement $x^{\perp}$ is given by $1-x$.

\item For $r\in [0,1]$ and $x\in [0,1]_{A}$ the scalar multiplications
  $rx$ and $(1-r)x$ are positive, and their sum is $x \leq 1$.  Hence
  $rx\leq 1$ and thus $rx\in[0,1]_{A}$.
\end{itemize}

\noindent Each PU-map of $C^*$-algebras $f\colon A \rightarrow B$
preserves $\leq$ and thus restricts to $[0,1]_{A} \rightarrow
[0,1]_{B}$. This restriction is a map of effect modules. Hence we get
a ``predicate'' functor $\CstarPU \rightarrow \EMod$.

\begin{lem}
\label{CstarPredFunLem}
The functor $[0,1]_{(-)} \colon \CstarPU \rightarrow \EMod$ is full
and faithful.
\end{lem}

\proof
Any PU-map $f\colon A \rightarrow B$ is completely determined (and
defined by) its action on $[0,1]_{A}$: for a non-zero positive element
$x\in A$ we use $x \leq \|x\|\, 1$ and thus $\frac{1}{\|x\|}\,x \in
[0,1]_{A}$ to see that $f(x) = \|x\|\,f(\frac{1}{\|x\|}\,x)$. An
arbitrary element $y\in A$ can be written uniquely as linear sum of four
positive elements (see Lemma~\ref{CstarMapLem}), determining
$f(y)$. 
\qed

The (finite, discrete probability) distribution monad $\Dst\colon
\Sets \rightarrow \Sets$ sends a set $X$ to the set $\Dst(X) =
\set{\varphi\colon X \rightarrow [0,1]}{\supp(\varphi) \mbox{ is
    finite, and }\sum_{x}\varphi(x) = 1}$, where $\supp(\varphi) =
\set{x}{\varphi(x)\neq 0}$.  Such an element $\varphi\in\Dst(X)$ may
be identified with a finite, formal convex sum $\sum_{i}r_{i}x_{i}$
with $x_{i}\in X$ and $r_{i}\in[0,1]$ satisfying
$\sum_{i}r_{i}=1$. The unit $\eta\colon X \rightarrow \Dst(X)$ and
multiplication $\mu \colon \Dst^{2}(X) \rightarrow \Dst(X)$ of this
monad are given by singleton/Dirac convex sum and by matrix
multiplication:
$$\begin{array}{rclcrcl}
\eta(x)
& = &
1x
& \qquad &
\mu(\Phi)(x)
& = &
\sum_{\varphi} \Phi(\varphi)\cdot \varphi(x).
\end{array}$$

\noindent A \emph{convex} set is an Eilenberg-Moore algebra of this
monad: it consists of a carrier set $X$ in which actual sums
$\sum_{i}r_{i}x_{i} \in X$ exist for all convex combinations. We write
$\Conv = \EM(\Dst)$ for the category of convex sets, with ``affine''
functions preserving convex sums.

Effect modules and convex sets are related via a basic
adjunction~\cite{JacobsM13a}, obtained by ``homming into $[0,1]$'', as
in:
\begin{equation}
\label{EModConvAdjDiag}
\xymatrix{
\op{\EMod}\ar@/^1.5ex/[rr]^-{\EMod(-,[0,1])} & \top &
   \Conv\ar@/^1.5ex/[ll]^-{\Conv(-,[0,1])}
}
\end{equation}


\section{Set-theoretic computations in \texorpdfstring{$C^*$}{C*}-algebras}\label{SetsSec}

For a set $X$, a function $f\colon X \rightarrow \C$ is called
\emph{bounded} if $|f(x)| \leq s$, for some $s\in\R_{\geq 0}$. We
write $\ell^{\infty}(X)$ for the set of such bounded functions. Notice
that if $X$ is finite, any function $X \rightarrow \C$ is bounded, so
that $\ell^{\infty}(X) = \C^{X}$.

Each $\ell^{\infty}(X)$ is a commutative $C^*$-algebra, with pointwise
addition, multiplication and involution, and with the uniform/supremum
norm:
$$\begin{array}{rcl}
\unifnorm{f}
& = &
\inf\setin{s}{\R_{\geq 0}}{\all{x}{|f(x)|\leq s}}.
\end{array}$$

\noindent In fact it is a typical example of a commutative
$W^*$-algebra, but we do not require this fact. This yields a functor
$\ell^{\infty} \colon \Sets \rightarrow \op{(\CCstarMIU)}$, where for
$h\colon X \rightarrow Y$ we have $\ell^{\infty}(h) = (-) \after h
\colon \ell^{\infty}(Y) \rightarrow \ell^{\infty}(X)$; it preserves
the (pointwise) operations. We have the following result.

\begin{prop}
\label{SetsCstarAdjProp}
The functor $\ell^{\infty} \colon \Sets \rightarrow
\op{(\CCstarMIU)}$ is left adjoint to the multiplicative states
functor $\MStat \colon \op{(\CCstarMIU)} \rightarrow \Sets$. In
combination with the adjunctions from Theorem~\ref{GelfandThm} we
get a situation:
$$\xymatrix@R+.5pc{
\quad\CH\;\ar@/^1.2ex/[dr]_{\dashv\!}\ar@/^1ex/[rr]^-{C} & \simeq & 
   \op{(\CCstarMIU)}\ar@/^1ex/[ll]^-{\MStat}
      \ar@/^1.2ex/[dl]_{\dashv\!}^{\MStat} \\
& \Sets\ar@/^1.5ex/[ul]^{\Ult}\ar@/^1.2ex/[ur]^(0.3){\ell^{\infty}\!\!} &
}$$

\noindent By composition and uniqueness of adjoints we get:
$$\begin{array}{rclcrcl}
C \after \Ult
& \cong &
\ell^{\infty}
& \qquad\mbox{and also}\qquad &
\MStat \after \ell^{\infty}
& \cong &
\Ult.
\end{array}$$
\end{prop}

\begin{proof}
Note that $\MStat$ is used in two different senses in the above
diagram, in one case with a compact Hausdorff topology, and in the
other case simply as a set. The adjunction involving $\ell^\infty$ and
$\MStat$ is for $\MStat$ as a set. We show this adjunction using the
universal property of the unit of an adjunction. We define the unit
$\eta_X \colon X \rightarrow \MStat(\ell^\infty(X))$, where $X \in
\Sets$, as
$$\begin{array}{rcl}
\eta_X(x)(a) & = & a(x),
\end{array}$$

\noindent where $a \in \ell^\infty(X)$. Then $\eta_X(x)$ is a
multiplicative state on $\ell^\infty(X)$ because the vector space
structure, multiplication and multiplicative unit are defined
pointwise. To show the naturality square for $\eta$ commutes, we must
show that for all $f \colon X \rightarrow Y$ in $\Sets$,
$\MStat(\ell^\infty(f)) \after \eta_X = \eta_Y \after f$. If we take $x
\in X$ and $b \in \ell^\infty(Y)$, we have:
$$\begin{array}{rcl}
\big(\MStat(\ell^\infty(f)) \after \eta_X\big)(x)(b)
& = &
\MStat(\ell^\infty(f))(\eta_X(x))(b) \\
& = &
(\eta_X(x) \after \ell^\infty(f))(b) \\
& = & 
\eta_X(x)(\ell^\infty(f)(b)) \\
& = &
\eta_X(x)(b \after f) \\
& = & 
b(f(x)) \\
& = &
\eta_Y(f(x))(b) \\
& = &
(\eta_Y \after f)(x)(b).
\end{array}$$

We now show this natural transformation satisfies the universal
property making it the unit of the adjunction. Let $X \in \Sets$, $B
\in \CCstarMIU$ and $f \colon X \rightarrow \MStat(B)$. Define $g: B
\rightarrow \ell^\infty(X)$ as $g(b)(x) = f(x)(b)$. We must show that
$g(b)$ is an element of $\ell^\infty(X)$, \emph{i.e.} that it is
bounded. For all $x \in X$, $f(x)$ is a multiplicative state, hence a
state, so by \cite[Proposition 2.1.4]{Dixmier77} we have $\|f(x)\| =
1$, and so $|g(b)(x) = |f(x)(b)| \leq \|f(x)\|\|b\| =
\|b\|$. Therefore $\|b\|$ is a bound for $g(b)$, showing that it is a
bounded function. The fact that $g$ is an MIU map is easily deduced
from the fact that $f(x)$ is a multiplicative state for all $x$ (it
would fail if $f(x)$ were only a state).

We must now show that
\[
\xymatrix@C+2pc{
X \ar[r]^-{\eta_X} \ar[rd]_f & \MStat(\ell^\infty(X)) \ar[d]^{\MStat(g)} \\
 & \MStat(B)
}
\]
commutes. Taking $x \in X$ and $b \in B$, we see
$$\begin{array}{rcl}
\MStat(g)(\eta_X(x))(b) 
& = &
(\eta_X(x) \after g)(b) \\
& = &
\eta_X(x)(g(b)) \\
& = &
g(b)(x) \\
& = &
f(x)(b),
\end{array}$$
and hence the unit diagram commutes.

To show the uniqueness of $g$, suppose there were $h \colon B
\rightarrow \ell^\infty(X)$ that also made the unit diagram
commute. By evaluating $\MStat(h)(\eta_X(x))(b)$ we would obtain
$g(b)(x) = h(b)(x)$. Since $g(b)$ and $h(b)$ are elements of
$\ell^\infty(X)$ and hence functions, this implies $g(b) = h(b)$ by
extensionality, and we can then conclude that $g = h$, as required. We
have now shown that $\ell^\infty$ is a left adjoint to $\MStat$. The
other two adjunctions are simply the Stone-\v{C}ech compactification
of a set and Gelfand duality (which is even an equivalence).

Since the triangle consisting of $\MStat$, in both forms, and the
forgetful functor $\CH \rightarrow \Sets$ commutes, the triangle for
$\ell^\infty, \Ult$ and $C$ commutes up to isomorphism, \emph{i.e.}
$\ell^\infty \cong C \after \Ult$ by uniqueness of adjoints.
\end{proof}

When we restrict to the full subcategory $\FinSets \hookrightarrow
\Sets$ of finite sets we obtain a functor $\ell^{\infty} = \C^{(-)}
\colon \FinSets \rightarrow \op{(\FdCCstarMIU)}$. The next result is
then a well-known special case of Gelfand duality
(Theorem~\ref{GelfandThm}). We elaborate the proof in some detail
because it is important to see where the preservation of
multiplication plays a role.

\begin{prop}
\label{FinSetsGelfandProp}
The functor $\C^{(-)} \colon \FinSets \rightarrow \op{(\FdCCstarMIU)}$
is an equivalence of categories.
\end{prop}

\proof
It is easy to see that the functor $\C^{(-)}$ is faithful. The
crucial part is to see that it is full. So assume we have two finite
sets, seen as natural numbers $n,m$, and a MIU-homomorphism $h\colon
\C^{m} \rightarrow \C^{n}$. For $j\in m$, let $\ket{j}\in\C^{m}$ be
the standard base vector with $1$ at the $j$-th position and $0$
elsewhere. Since this $\ket{j}$ is positive, so is $h(\ket{j})$, and
thus we may write it as $h(\ket{j}) = (r_{1j}, \ldots, r_{nj})$, with
$r_{ij} \in \R_{\geq 0}$. Because $\ket{j} \cdot \ket{j} = \ket{j}$,
and $h$ preserves multiplication, we get $h(\ket{j}) \cdot h(\ket{j})
= h(\ket{j})$, and thus $r_{ij}^{2} = r_{ij}$. This means $r_{ij}\in
\{0,1\}$, so that $h$ is a (binary) Boolean matrix. But $h$ is also
unital, and so:
\begin{equation}
\label{UnitalBasisEqn}
\begin{array}{rcccccl}
1
& = &
h(1)
& = &
h(\ket{1}+\cdots+\ket{m})
& = & 
h(\ket{1})+ \cdots + h(\ket{m}).
\end{array}
\end{equation}

\noindent For each $i\in n$ there is thus precisely one $j\in m$ with
$r_{ij}=1$ --- so that $h$ is a ``functional'' Boolean matrix. This
yields the required function $f\colon n \rightarrow m$ with $\C^{f} =
h$.

Corollary~\ref{FinGelfandCor} says that the functor $\C^{(-)} \colon
\FinSets \rightarrow \op{(\FdCCstarMIU)}$ is essentially surjective on
objects, and thus an equivalence. 
\qed

This proof demonstrates that preservation of multiplication, as
required for ``MIU'' maps, is a rather strong condition. We make this
more explicit.

\begin{cor}
\label{MIUStateCor}
For $n\in\NNO$ we have $\MStat(\C^{n}) \cong n$.
\end{cor}

\proof
By identifying $n\in\NNO$ with the $n$-element set $n =
\{0,1,\ldots,n-1\} \in\FinSets$, we get by
Proposition~\ref{FinSetsGelfandProp}, $\MStat(\C^{n}) = \HomMIU(\C^{n},
\C) \cong \FinSets(1, n) \cong n$. 
\qed

\auxproof{
\begin{rem}
\label{NondetRem}
Looking at the proof of Proposition~\ref{FinSetsGelfandProp} one sees
that the ``non-deterministic'' computations $n \rightarrow \Pow(m)$
correspond to linear maps $\C^{m} \rightarrow \C^{n}$ that are
positive and preserve squares of the standard base vectors
$\ket{j}$. It is unclear if there is also a more algebraic,
non-basis-dependent, formulation of this property.

More categorically, let us write $\U$ for the ultrafilter monad on
$\Sets$. The category $\CH$ of compact Hausdorff spaces is the
category of Eilenberg-Moore algebras of $\mathcal{U}$. Starting from
the Kleisli category $\Kl(\U)$ of this monad we then get a full and
faithful functor:
$$\xymatrix{
\Kl(\U)\ar[rr]^-{\textrm{full \& faithful}} & &
   \EM(\U) = \CH\ar[rr]^-{\Cont(-,\C)}_-{\simeq} & & 
   \op{(\CCstarMIU)}
}$$

\noindent If we now use that $\U(X)$ for a finite set $X$ is just
powerset $\Pow(X)$, then we have another way of embedding
non-deterministic computations (between finite sets) in
$C^*$-algebras, namely via $n\mapsto \C^{\Pow(n)}$.\qed
\end{rem}
}

\section{Discrete probabilistic computations in
  \texorpdfstring{$C^*$}{C*}-algebras}\label{DiscProbSec}

We turn to probabilistic computations and will see that we remain in
the world of commutative $C^{*}$-algebras, but with PU-maps (positive
unital) instead of MIU-maps. Recall that the set of states $\Stat(A)$
of a $C^*$-algebra $A$ contains the PU-maps $A \rightarrow \C$.

We summarize here the definition of the expectation monad given in
\cite{JacobsM12b}. If $[0,1]^X$ is the effect module of functions from
$X$ to $[0,1]$ with pointwise operations, $\E(X) = \EMod([0,1]^X,
     [0,1])$. The unit $\eta_X \colon X \rightarrow \E(X)$ is
     evaluation, defined as $\eta_X(x)(f) = f(x)$ for $f \in
     [0,1]^X$. The multiplication $\mu_X \colon \E^2(X) \rightarrow
     \E(X)$ is defined for $h \in [0,1]^{\E(X)} \rightarrow [0,1]$, $p
     \in [0,1]^X$ as
$$\begin{array}{rcl}
\mu_X(h)(p)
& = &
h\big(\lamin{k}{\E(X)}{k(p)}\big).
\end{array}$$

\begin{lem}
\label{ExpectationLem}
Sending a set $X$ to the set of states of the $C^{*}$-algebra
$\ell^{\infty}(X)$ yields the (underlying functor of the) expectation monad $\E$
from~\cite{JacobsM12b}: the mapping $X \mapsto
\Stat(\ell^{\infty}(X))$ is isomorphic to the expectation monad $\E
\colon \Sets \rightarrow \Sets$, defined in~\cite{JacobsM12b} via
effect module homomorphisms: $\E(X) = \EMod\big([0,1]^{X},
[0,1]\big)$.

As a result, $\Stat(\C^{n}) \cong \Dst(n)$, for $n\in\NNO$, where
$\Dst(n)$ is the standard $n$-simplex.
\end{lem}


\proof
The predicate/effect functor $[0,1]_{(-)} \colon \CstarPU \rightarrow
\EMod$ is full and faithful by Lemma~\ref{CstarPredFunLem}, and so:
$$\begin{array}{rcl}
\Stat(\ell^{\infty}(X))
\hspace*{\arraycolsep} = \hspace*{\arraycolsep}
\HomPU\big(\ell^{\infty}(X), \C\big)
& \cong &
\EMod\big([0,1]_{\ell^{\infty}(X)}, [0,1]_{\C}\big) \\
& = &
\EMod\big([0,1]^{X}, [0,1]\big)
\hspace*{\arraycolsep} = \hspace*{\arraycolsep}
\E(X).
\end{array}$$

\noindent The isomorphism $\smash{\alpha \colon \HomPU(\C^{n}, \C)
  \conglongrightarrow \Dst(n)}$ follows because the expectation and
distribution monad coincide on finite sets,
see~\cite{JacobsM12b}. Explicitly, it is given by $\alpha(h) =
\lamin{i}{n}{h(\ket{i})}$ and $\alpha^{-1}(\varphi)(v) =
\sum_{i}\varphi(i) \cdot v(i)$. \qed

\auxproof{
We check that $\alpha$ is well-defined: since $\ket{i}\in \C^{n}$ is
positive, we get $h(\ket{i}) \geq 0$; and:
$$\begin{array}{rcl}
\sum_{i}\alpha(h)(i)
& = &
\sum_{i} h(\ket{i}) \\
& = &
h(\sum_{i}\ket{i}) \\
& = &
h(1) \\
& = &
1.
\end{array}$$

In the reverse direction we have to show that $\alpha^{-1}(\varphi)
\colon \C^{n} \rightarrow \C$ is a PU-map:
$$\begin{array}{rcl}
\alpha^{-1}(\varphi)(1)
& = &
\sum_{i}\varphi(i)\cdot 1(i) \\
& = &
\sum_{i}\varphi(i) \\
& = &
1 \\
\alpha^{-1}(\varphi)(v^{*}v)
& = &
\sum_{i}\varphi(i)\cdot (v^{*}v)(i) \\
& \geq &
0 \qquad \mbox{as sum of positive numbers.}
\end{array}$$

Finally:
$$\begin{array}{rcl}
\big(\alpha \after \alpha^{-1}\big)(\varphi)(i)
& = &
\alpha\big(\alpha^{-1}(\varphi)\big)(i) \\
& = &
\alpha^{-1}(\varphi)(\ket{i}) \\
& = &
\sum_{j} \varphi(j)\cdot \ket{i}(j) \\
& = &
\varphi(i) \\
\big(\alpha^{-1} \after \alpha\big)(h)(v)
& = &
\alpha^{-1}\big(\alpha(h)\big)(v) \\
& = &
\sum_{i} \alpha(h)(i) \cdot v(i) \\
& = &
\sum_{i} h(\ket{i}) \cdot v(i) \\
& = &
\sum_{i} h(v(i) \cdot \ket{i}) \\
& = &
h(\sum_{i} v(i) \cdot \ket{i}) \\
& = &
h(v).
\end{array}$$
}

The unit $\eta$ and multiplication $\mu$ structure on $\E(X) \cong
\HomPU(\ell^{\infty}(X), \C)$ is very much like for ``continuation''
or ``double dual'' monads, see~\cite{Kock70b,Moggi91a,Jacobs12}, with:
$$\hspace*{-.5em}\xymatrix@R-2pc@C-.8pc{
X\ar[r]^-{\eta} & \HomPU(\ell^{\infty}(X), \C)
\!\! & \!\!
\HomPU\Big(\ell^{\infty}\big(\HomPU(\C^{X}, \C)\big), \C\Big)\ar[r]^-{\mu} & 
   \HomPU(\ell^{\infty}(X), \C) \\
x\ar@{|->}[r] & \lam{v}{v(x)}
&
g\ar@{|->}[r] & \lam{v}{g\big(\lam{h}{h(v)}\big)}.
}$$

For an arbitrary monad $T = (T, \eta, \mu)$ on a category $\cat{B}$ we
write $\Kl(T)$ for the Kleisli category of $T$. Its objects are the
same as those of $\cat{B}$, but its maps $X \rightarrow Y$ are the
maps $X \rightarrow T(Y)$ in $\cat{B}$. The unit $\eta\colon X
\rightarrow T(X)$ is the identity map $X \rightarrow X$ in $\Kl(T)$;
and composition of $f\colon X \rightarrow Y$ and $g\colon Y
\rightarrow Z$ in $\Kl(T)$ is given by $g \klafter f = \mu \after T(g)
\after f$. Maps in such a Kleisli category are understood as
computations with outcomes of type $T$, see~\cite{Moggi91a}.  For a
monad $T\colon \Sets \rightarrow \Sets$ we write $\KlN(T)
\hookrightarrow \Kl(T)$ for the full subcategory with numbers
$n\in\NNO$ as objects, considered as $n$-element sets.

\begin{prop}
\label{KlEToStarProp}
The expectation monad $\E(X) \cong \HomPU(\ell^{\infty}(X), \C)$ gives
rise to a full and faithful functor:
\begin{equation}
\label{KlEToStarFun}
\vcenter{\xymatrix@R-2pc{
\Kl(\E)\ar[rr]^-{\starfun{\E}} & & \op{(\CCstarPU)} \\
X\ar@{|->}[rr] & & \ell^{\infty}(X) \\
\big(X\stackrel{f}{\rightarrow} \E(Y)\big)\ar@{|->}[rr] & &
   \lamin{v}{\ell^{\infty}(Y)}{\lamin{x}{X}{f(x)(v)}}.
}}
\end{equation}
\end{prop}

\proof
First we need to see that $\starfun{\E}(f)$ is well-defined: the
function $\starfun{\E}(f)(v) \colon X \rightarrow \C$ must be
bounded. We can apply Lemma~\ref{CstarMapLem} to the function $f(x)
\in \HomPU(\ell^{\infty}(Y), \C)$; it yields $\|f(x)(v)\| \leq
4\|v\|$. This holds for each $x\in X$, so that
$|\starfun{\E}(f)(v)(x)| = |f(x)(v)|$ is bounded by $4\|v\|$.  Next,
the map $\starfun{\E}(f)$ is a PU-map of $C^*$-algebras via the
pointwise definitions of the relevant constructions.

We check that $\starfun{\E}$ preserves (Kleisli) identities and
composition:
$$\begin{array}{rcl}
\starfun{\E}(\idmap)(v)(x)
& = &
\starfun{\E}(\eta)(v)(x) \\
& = &
\eta(x)(v) \\
& = &
v(x) \\
\starfun{\E}(g \klafter f)(v)(x)
& = &
(g \klafter f)(x)(v) \\
& = &
\mu\Big(\E(g)(f(x))\Big)(v) \\
& = &
\E(g)(f(x))\big(\lam{w}{w(v)}\big) \\
& = &
f(x)\big((\lam{w}{w(v)}) \after g\big) \\
& = &
f(x)\big(\lam{y}{g(y)(v)}\big) \\
& = &
f(x)\big(\starfun{\E}(g)(v)\big) \\
& = &
\starfun{\E}(f)\big(\starfun{\E}(g)(v)\big)(x) \\
& = &
\big(\starfun{\E}(f) \after \starfun{\E}(g)\big)(v)(x).
\end{array}$$

\noindent Further, $\starfun{\E}$ is obviously faithful, and it is
full since for $h\colon \ell^{\infty}(Y) \rightarrow \ell^{\infty}(X)$
in $\CCstarPU$ we can define $f\colon X \rightarrow
\HomPU(\ell^{\infty}(Y), \C)$ by $f(x)(v) = h(v)(x)$.  Then each
$f(x)$ is a PU-map of $C^*$-algebras. \qed

\auxproof{
Preservation of the pointwise operations:
$$\begin{array}{rcl}
\starfun{\E}(f)(v+v')(x)
& = &
f(x)(v+v') \\
& = &
f(x)(v) + f(x)(v') \\
& = &
\starfun{\E}(f)(v)(x) + \starfun{\E}(f)(v')(x).
\end{array}$$

For $h\colon \C^{Y} \rightarrow \C^{X}$ in $\CCstarPU$ we define
$f\colon X \rightarrow \HomPU(\C^{Y}, \C)$ by $f(x)(v) = h(v)(x)$
and get for each $f(x)$ a PU-map of $C^*$-algebras, since:
$$\begin{array}{rcl}
f(x)(1)
& = &
h(1)(x) \\
& = &
1(x) \\
& = & 
x \\
f(x)(v^{*}v)
& = &
h(v^{*}v)(x) \\
& \geq &
0 \qquad \mbox{since } h(v^{*}v) \geq 0 \mbox{ in } \C^{X}.
\end{array}$$
}

We turn to the finite case, like in the previous section. We do so by
considering the Kleisli category $\KlN(\E)$ obtained by restricting to
objects $n\in\NNO$. Since the expectation monad $\E$ and the
distribution monad $\Dst$ coincide on finite sets, we have $\KlN(\E)
\cong \KlN(\Dst)$. Maps $n \rightarrow m$ in this category are
probabilistic transition matrices $n \rightarrow \Dst(m)$. This category has been investigated also in \cite{Fritz09}.  The
following equivalence is known, see \textit{e.g.}~\cite{Maassen10},
although possibly not in this categorical form.

\begin{prop}
\label{FinKlDToStarFunProp}
The functor $\starfun{\E}$ from~\eqref{KlEToStarFun} restricts in the
finite case to an equivalence of categories:
\begin{equation}
\label{FinKlDToStarFun}
\vcenter{\xymatrix{
\KlN(\Dst)\ar[rr]^-{\starfun{\Dst}}_-{\simeq} & & \op{(\FdCCstarPU)}
}}
\end{equation}

\noindent It is given by $\starfun{\Dst}(n) = \C^{n}$ and
$\smash{\starfun{\Dst}\big(n \stackrel{f}{\rightarrow} \Dst(m)\big)
= \lamin{v}{\C^{m}}{\lamin{i}{n}{\sum\limits_{j\in m} f(i)(j)\cdot v(j)}}}$.
\end{prop}

This equivalence~\eqref{FinKlDToStarFun} may be read as: the category
$\FdCCstarPU$ of finite-dimensional commutative $C^*$-algebras, with
positive unital maps, is the \emph{Lawvere theory} of the distribution
monad $\Dst$.

\proof
Fullness and faithfulness of the functor $\starfun{\Dst}$ follow from
Proposition~\ref{KlEToStarProp}, using the isomorphism $\HomPU(\C^{n},
\C) \cong \Dst(n)$ from Lemma~\ref{ExpectationLem}.  This functor
$\starfun{\Dst}$ is essentially surjective on objects by
Corollary~\ref{FinGelfandCor}, using the fact that a MIU-map is a PU-map. \qed

\auxproof{
Explicit, older proof:

Assume two functions $f,g\colon n \rightarrow \Dst(m)$ with $P(f) =
P(g) \colon \C^{m} \rightarrow \C^{n}$. The standard base vector
$\ket{j}\in\C^{m}$, for $j\in m$, yields, for $i\in n$:
$$\begin{array}{rcccccl}
f(i)(j)
& = &
\bra{i}P(f)\ket{j}
& = &
\bra{i}P(g)\ket{j}
& = &
g(i)(j).
\end{array}$$

\noindent Hence $f=g$, so that $P$ is a faithful functor. It is also
full, since a $C^*$-algebra map $h\colon \C^{m} \rightarrow \C^{n}$
consists of an $n\times m$ matrix $M = (M_{ij})$ preserving unit and
positivity.  Since each $\ket{j}\in\C^{m}$, for $j\in m$, is positive,
so is the vector $M\ket{j} = (M_{1j}, \ldots, M_{nj})\in\C^{n}$. This
means each $M_{ij}\in\C$ is positive. Since $M$ is unital we get
$M_{i1} + \cdots + M_{im} = 1$ for each $i\in n$,
via~\eqref{UnitalBasisEqn}. Thus we can define a Kleisli map $f\colon
n \rightarrow \Dst(m)$ by $f(i)(j) = M_{ij}\in [0,1]$, with
$\sum_{j}f(i)(j) = 1$. Then $P(f) = M$ since:
$$\begin{array}{rcccccccl}
P(f)(v)(i)
& = &
\sum_{j} f(i)(j)\cdot v(j) 
& = &
\sum_{j} M_{ij}\cdot v(j) 
& = &
(M v)(i)
& = &
h(v)(i).
\end{array}$$

We briefly check that $\starfun{\Dst}$ is a functor:
$$\begin{array}{rcl}
\starfun{\Dst}(\idmap)(v)(x)
& = &
\sum_{y}\eta(x)(y)\cdot v(y) \\
& = &
v(x) \\
\starfun{\Dst}(g \klafter f)(v)(x)
& = &
\sum_{z} (g \klafter f)(x)(z)\cdot v(z) \\
& = &
\sum_{z} \big(\sum_{y} f(x)(y)\cdot g(y)(z)\big) \cdot v(z) \\
& = &
\sum_{y,z} f(x)(y)\cdot g(y)(z) \cdot v(z) \\
& = &
\sum_{y} f(x)(y)\cdot \big(\sum_{z} g(y)(z) \cdot v(z)\big) \\
& = &
\sum_{y} f(x)(y)\cdot \starfun{\Dst}(g)(v)(y) \\
& = &
\starfun{\Dst}(f)\big(\starfun{\Dst}(g)(v)\big)(x) \\
& = &
\big(\starfun{\Dst}(f) \after \starfun{\Dst}(g)\big)(v)(x).
\end{array}$$
}

\section{Continuous probabilistic computations}\label{ContProbSec}

The question arises if the full and faithful functor $\Kl(\E)
\rightarrow \op{(\CCstarPU)}$ from Proposition~\ref{KlEToStarProp} can
be turned into an equivalence of categories, but not just for the
finite case like in Proposition~\ref{FinKlDToStarFunProp}. In order to
make this work we have to lift the expectation monad $\E$ on $\Sets$
to the category $\CH$ of compact Hausdorff spaces. As lifting we use
what we call the \emph{Radon} monad $\Rdn$, defined on $X\in\CH$ as:
\begin{equation}
\label{RdnEqn}
\begin{array}{rcccl}
\Rdn(X)
& = &
\Stat(C(X))
& = &
\HomPU\big(C(X), \, \C\big),
\end{array}
\end{equation}

\noindent where, as usual, $C(X) = \set{f\colon X \rightarrow \C}{f
  \mbox{ is continuous}}$; notice that the functions $f\in C(X)$ are
automatically bounded, since $X$ is compact. We have implicitly
applied the forgetful functor from $\CCL \rightarrow \CH$ to make
$\Rdn$ into an endofunctor of $\CH$.  The elements of $\Rdn(X)$ are
related to measures in the following way. If $\mu$ is a probability
measure on the Borel sets of $X$, integration of continuous functions
with respect to $\mu$ gives a function $\int_X - d\mu \in \Rdn(X)$. A
Radon probability measure, or an inner regular probability measure, is
one such that $\mu(S) = \sup_{K \subseteq S}\mu(K)$ where $K$ ranges
over compact sets. The map from measures to elements of $\Rdn(X)$ is a
bijection~\cite[Thm.~2.14]{Rudin87}, and accordingly we shall
sometimes refer to elements of $\Rdn(X)$ as measures. Therefore the
Radon monad can be considered as a variant of the Giry monad. In fact
there are two Giry monads, one on measurable spaces and one on Polish
spaces. The Radon monad differs from the Giry monad on measurable
spaces in that it uses the topology of a space, and that in the case
of a space that is not a standard Borel space there can be non-Radon
measures \cite[434K (d), page 192]{fremlin} \cite[\S 53.10, page
  231]{Halmos50}. The Radon monad differs from the Giry monad on
Polish spaces essentially only in the choice of spaces, and on compact
Polish spaces they agree, as the topology Giry used is the same as the
weak-* topology, and Polish spaces do not admit any non-Radon Borel
probability measures.\cite[Theorems 1.1 and 1.4]{Billingsley68}.

This Radon monad $\Rdn$ is not new: we shall see later that it occurs
in~\cite[Theorem 3]{Swirszcz74} as the monad of an adjunction
(``probability measure'' is used to mean ``Radon probability measure''
in that article). It has been used more recently in
\cite{Mislove2012}. However, our duality result below ---
Theorem~\ref{KlRToStarThm} --- is not known in the literature.

From Proposition~\ref{StatesProp} it is immediate that $\Rdn(X)$ is
again a compact Hausdorff space. The unit $\eta\colon X \rightarrow
\Rdn(X)$ and multiplication $\mu \colon \Rdn^{2}(X) \rightarrow
\Rdn(X)$ are defined as for the expectation monad, namely as
$\eta(x)(v) = v(x)$ and $\mu(g)(v) = g\big(\lam{h}{h(v)}\big)$.  We
check that $\eta$ is continuous. Recall from the proof of
Proposition~\ref{StatesProp} that a basic open in $\Rdn(X)$ is of the
form $\ev_{s}^{-1}(U) = \setin{h}{\Rdn(X)}{h(s)\in U}$, where
$s\in C(X)$ and $U\subseteq\C$ is open. Then:
$$\begin{array}{rcccccl}
\eta^{-1}\big(\ev_{s}^{-1}(U)\big)
& = &
\setin{x}{X}{\eta(x)(s)\in U}
& = &
\setin{x}{X}{s(x)\in U}
& = &
s^{-1}(U).
\end{array}$$

\noindent The latter is an open subset of $X$ since $s\colon X
\rightarrow\C$ is a continuous function.

\auxproof{
Also the multiplication is continuous: 
$$\begin{array}{rcl}
\mu^{-1}\big(\ev_{s}^{-1}(U)\big)
& = &
\setin{g}{\Rdn^{2}(X)}{\mu(g)(s) \in U} \\
& = &
\setin{g}{\Rdn^{2}(X)}{g(\lam{h}{h(s)}) \in U} \\
& = &
\setin{g}{\Rdn^{2}(X)}{\ev_{\lam{h}{h(s)}}(g) \in U} \\
& = &
\ev_{\lam{h}{h(s)}}^{-1}(U).
\end{array}$$
}

We are now ready to state our main, new duality result. It may be
understood as a probabilistic version of Gelfand duality, for
commutative $C^*$-algebras with PU maps instead of the MIU maps
originally used (see Theorem \ref{GelfandThm}).

\begin{thm}
\label{KlRToStarThm}
The Radon monad~\eqref{RdnEqn} yields an equivalence of
categories:
$$\begin{array}{rcl}
\Kl(\Rdn) 
& \simeq &
\op{(\CCstarPU)}.
\end{array}$$
\end{thm}

\proof
We define a functor $\starfun{\Rdn} \colon \Kl(\Rdn) \rightarrow
\op{(\CCstarPU)}$ like in~\eqref{KlEToStarFun}, namely by:
$$\begin{array}{rclcrcl}
\starfun{\Rdn}(X)
& = &
C(X) 
& \qquad &
\starfun{\Rdn}(f)
& = &
\lam{v}{\lam{x}{f(x)(v)}}.
\end{array}$$

\noindent Since $f\colon X \rightarrow \Rdn(Y)$ is itself continuous,
so is $f(-)(v) \colon X \rightarrow \C$.

\auxproof{
Let $f\colon X \rightarrow \Rdn(Y)$. Then for and open set $U\subseteq \C$
we obtain an open set:
$$\begin{array}{rcl}
\big(f(-)(v)\big)^{-1}(U)
& = &
\setin{x}{X}{f(x)(v) \in U} \\
& = &
\setin{x}{X}{\ev_{v}(f(x)) \in U} \\
& = &
\setin{x}{X}{f(x) \in \ev_{v}^{-1}(U)} \\
& = &
f^{-1}\big(\ev_{v}^{-1}(U)\big)
\end{array}$$
}

The fact that $\starfun{\Rdn}$ is a full and faithful functor follows
as in the proof of Proposition~\ref{KlEToStarProp}. This functor is
essentially surjective on objects by ordinary Gelfand duality
(Theorem~\ref{GelfandThm}). 
\qed

We investigate the Radon monad $\Rdn$ a bit further, in particular its
relation to the distribution monad $\Dst$ on $\Sets$.

\begin{lem}
\label{MonadFunctorLemma}
There is a map of monads $(U,\tau) \colon \Rdn \rightarrow \Dst$ in:
$$\xymatrix{
\CH\ar@(l,u)^{\Rdn}\ar[rr]^-{U} & & \Sets\ar@(u,r)^{\Dst}
& & \Dst U \ar@{=>}[r]^-{\tau} & U\Rdn
}$$

\noindent where $U$ is the forgetful functor and $\tau$ commutes
appropriately with the units and multiplications of the monads $\Dst$
and $\Rdn$. (Such a map is called a ``monad functor''
in~\cite[\S 1]{Street72}.)

As a result the forgetful functor lifts to the associated categories
of Eilenberg-Moore algebras:
$$\xymatrix@R-2pc{
\EM(\Rdn)\ar[rr] & & \EM(\Dst)\rlap{$\;=\Conv$} \\
\big(\Rdn(X)\stackrel{\alpha}{\rightarrow} X\big)\ar@{|->}[rr] & &
   \big(\Dst(UX)\stackrel{\tau}{\rightarrow} U\Rdn(X)
   \stackrel{U\alpha}{\rightarrow} UX\big)
}$$

\noindent Hence the carrier of an $\Rdn$-algebra is a convex compact
Hausdorff space, and every algebra map is an affine function. 
\end{lem}

\proof
For $X\in\CH$ and $\varphi\in\Dst(UX)$, that is for $\varphi\colon UX
\rightarrow [0,1]$ with finite support and $\sum_{x}\varphi(x)=1$, we
define $\tau(\varphi) \in U\Rdn(X)$ on $h\in C(X)$ as:
\begin{equation}
\label{DstToRdnEqn}
\begin{array}{rcl}
\tau(\varphi)(h)
& = &
\sum_{x} \varphi(x) \cdot h(x) \;\in\;\C.
\end{array}
\end{equation}

\noindent It is easy to see that $\tau$ is a linear map $C(X)
\rightarrow \C$ that preserves positive elements and the unit. Moreover,
it commutes appropriately with the units and multiplications. For
instance:
\[\big(\tau_{X} \after \eta^{\Dst}_{UX}\big)(x)(h)
 = 
\tau_{X}(1x)(h) 
 = 
h(x) 
 = 
U(\eta^{\Rdn}_{X})(x)(h).\eqno{\qEd}
\]\smallskip

\auxproof{
$$\begin{array}{rcl}
\big(U(\mu_{X}) \after \tau_{U\Rdn(X)} \after \Dst(\tau_{UX})\big)(\Phi)(h)
& = &
\mu_{X}\Big(\tau_{U\Rdn(X)}\big(\Dst(\tau_{UX})(\Phi)\big)\Big)(h) \\
& = &
\tau_{U\Rdn(X)}\big(\Dst(\tau_{UX})(\Phi)\big)\big(\lam{k}{k(h)}\big) \\
& = &
\sum_{k} \Dst(\tau_{UX})(\Phi)(k) \cdot k(h) \\
& = &
\sum_{k} \big(\sum_{\varphi\in\tau^{-1}(k)} \Phi(\varphi)\big) \cdot k(h) \\
& = &
\sum_{k, \varphi\in\tau^{-1}(k)} \Phi(\varphi) \cdot k(h) \\
& = &
\sum_{\varphi}  \Phi(\varphi) \cdot \tau(\varphi)(h) \\
& = &
\sum_{\varphi}  \Phi(\varphi) \cdot \big(\sum_{x}\varphi(x) \cdot h(x)\big) \\
& = &
\sum_{x,\varphi} \Phi(\varphi) \cdot \varphi(x) \cdot h(x) \\
& = &
\sum_{x} \big(\sum_{\varphi} \Phi(\varphi) \cdot \varphi(x)\big) \cdot h(x) \\
& = &
\sum_{x} \mu_{UX}(\Phi)(x)\cdot h(x) \\
& = &
\tau_{X}\big(\mu_{UX}(\Phi)\big)(h) \\
& = &
\big(\tau_{X} \after \mu_{UX}\big)(\Phi)(h).
\end{array}$$
}

\noindent The continuous dual space of $C(X)$ can be ordered using \eqref{ConePosetDef}, by taking the positive cone to be those linear functionals that map positive functions to positive numbers.

\begin{defi}
A state $\phi\in\Rdn(X)=\HomPU(C(X),\C)$ is a \emph{pure state} if for for each positive linear functional such that $\psi \leq \phi$, \emph{i.e.} such that $\phi - \psi$ is positive, there exists an $\alpha \in [0,1]$ such that $\psi = \alpha\phi$. \qed
\end{defi}

\begin{lem}
\label{ExtremeLemma}
For a compact Hausdorff space $X$, the subset of unit (or Dirac)
measures $\set{\eta(x)}{x\in X} \subseteq \Rdn(X)$ are pure states and hence is the set of
extreme points of the set of Radon measures $\Rdn(X)$ --- where
$\eta(x) = \eta^{\Rdn}(x) = \ev_{x} = \lam{h}{h(x)}$ is the unit of
the monad $\Rdn$.
\end{lem}

\proof
We rely on the basic fact, see~\cite[2.5.2, page 43]{Dixmier77}, that a measure is a Dirac measure iff it is a ``pure'' state. We prove
the above lemma by showing that the pure states are precisely the
extreme points of the convex set $\Rdn(X)$.


\begin{itemize}
\item If $\phi\in\Rdn(X)$ is a pure state, suppose $\phi = \alpha_1
  \phi_1 + \alpha_2 \phi_2$, a convex combination of two states
  $\phi_{i}\in \Rdn(X)$ with $\alpha_{i}\in[0,1]$ satisfying
  $\alpha_{1}+\alpha_{2} = 1$, where no two elements of $\{\phi,
  \phi_1, \phi_2 \}$ are the same. Then $\phi \geq \alpha_1 \phi_1$,
  since for a positive function $f\in C(X)$ one has $(\phi -
  \alpha_{1}\phi_{1})(f) = \alpha_{2}\phi_{2}(f) \geq 0$. Thus
  $\alpha_{1}\phi_{1} = \alpha\phi$, for some $\alpha\in[0,1]$, since
  $\phi$ is pure. Then $\alpha_{1} = \alpha_{1}\phi_{1}(1) =
  \alpha\phi(1) = \alpha$. If $\alpha_{1} = 0$, then $\alpha_{2} = 1$
  and so $\phi = \phi_{2}$.  If $\alpha_{1} > 0$, then $\phi =
  \phi_{1}$. Hence $\phi$ is an extreme point.


\item Suppose $\phi$ is an extreme point of $\Rdn(X)$,
  \textit{i.e.}~that $\phi = \alpha_1 \phi_1 + \alpha_2 \phi_2$
  implies $\phi_1$ or $\phi_2 = \phi$. Then if there is a positive
  linear functional $\psi \leq \phi$, we may take $\alpha_1 = \psi(1)
  \geq 0$; since $\alpha_{1} = \psi(1) \leq \phi(1) = 1$, we get
  $\alpha_{1} \in [0,1]$. If $\alpha_1 = 0$, then since $\|\psi\| =
  \psi(1) = 0$ we get $\psi = 0$ and $\psi = 0 \cdot \phi$.  If
  $\alpha_1 = 1$, then $(\phi - \psi)(1) = 0$, which since $\phi -
  \psi$ was assumed to be positive implies $\phi - \psi = 0$ and hence
  $\psi = 1 \cdot \phi$. Having dealt with those cases, we have that
  $\alpha_1 \in (0,1)$, and so we have a state $\phi_1 =
  \frac{1}{\alpha_1}\psi$. We may take $\alpha_2 = 1 - \alpha_{1} \in
  (0,1)$ and obtain a second state $\phi_2 = \frac{1}{\alpha_2}(\phi -
  \psi)$. By construction we have a convex decomposition of $\phi =
  \alpha_1\phi_1 + \alpha_2\phi_2$. Therefore either $\phi = \phi_{1}
  = \frac{1}{\alpha_1}\psi$ or $\phi = \phi_{2} =
  \frac{1}{\alpha_2}(\phi - \psi)$. In the first case, $\psi =
  \alpha_1 \phi$, making $\phi$ pure. But also in the second case
  $\phi$ is pure, since we have $\alpha_2 \phi = \phi - \psi$ and thus
  $\psi = (1- \alpha_2)\phi$. \qed
\end{itemize}

\begin{lem}
\label{DensityLemma}
Let $X$ be a compact Hausdorff space.
\begin{enumerate}
\item The maps $\tau_{X} \colon \Dst(UX) \rightarrow U\Rdn(X)$
  from~\eqref{DstToRdnEqn} are injective; as a result, the unit/Dirac
  maps $\eta\colon X \rightarrow \Rdn(X)$ are also injective.

\item The maps $\tau_{X} \colon \Dst(UX) \rightarrowtail U\Rdn(X)$ are
  dense.
\end{enumerate}
\end{lem}

\proof
For the first point, assume $\varphi,\psi\in\Dst(UX)$ satisfying
$\tau(\varphi) = \tau(\psi)$. We first show that the finite support
sets are equal: $\supp(\varphi) = \supp(\psi)$. Since $X$ is
Hausdorff, singletons are closed, and hence finite subsets too.
Suppose $\supp(\varphi) \not\subseteq \supp(\psi)$, so that $S =
\supp(\varphi) - \supp(\psi)$ is non-empty. Since $S$ and
$\supp(\psi)$ are disjoint closed subsets, there is by Urysohn's lemma
a continuous function $f\colon X \rightarrow [0,1]$ with $f(x)=1$ for
$x\in S$ and $f(x) = 0$ for $x\in \supp(\psi)$. But then
$\tau(\psi)(f) = 0$, whereas $\tau(\varphi)(f) \neq 0$.

Now that we know $\supp(\varphi) = \supp(\psi)$, assume $\varphi(x)
\neq \psi(x)$, for some $x\in\supp(\varphi)$. The closed subsets
$\{x\}$ and $\supp(\varphi) - \{x\}$ are disjoint, so there is, again
by Urysohn's lemma a continuous function $f\colon X \rightarrow [0,1]$
with $f(x) = 1$ and $f(y) = 0$ for all $y\in\supp(\varphi)$. But then
$\varphi(x) = \tau(\varphi)(f) = \tau(\psi)(f) = \psi(x)$,
contradicting the assumption.

We can conclude that the unit $X \rightarrow \Rdn(X)$ is also
injective, since its underlying function can be written as the composite
$U(\eta^\Rdn) = \tau \after \eta^\Dst \colon UX \rightarrowtail \Dst(UX)
\rightarrowtail U\Rdn(X)$, because $\tau$ is a map of monads.

To show that the image of $\tau_X$ is dense, we proceed as follows. By
Lemmas~\ref{ExtremeLemma} and~\ref{MonadFunctorLemma}, the extreme
points of $\Rdn(X)$ are
$$\begin{array}{rcl}
\set{\eta^{\Rdn}(x)}{x\in X}
& = &
\set{\tau\big(\eta^{\Dst}(x))}{x\in X}
\end{array}$$

\noindent and are thus in the image of $\tau\colon \Dst(UX)
\rightarrowtail U\Rdn(X)$. Since every convex combination of
$\eta^{\Rdn}(x)$ comes from a formal convex sum $\varphi \in
\Dst(UX)$, all convex combinations of extreme points are in the image
of $\tau_X$. Using Proposition \ref{StatesProp}, $\Rdn(X)$ can be considered an object of $\CCL$, i.e. a compact convex subset of a locally convex space. Accordingly, we may apply the Krein-Milman theorem~\cite[Proposition 7.4, page
  142]{Conway90} to conclude the set of convex combinations of extreme
points is dense. 
\qed

\begin{lem}
\label{AffineLem}
Let $X$ and $Y$ be compact Hausdorff spaces.  Each Eilenberg-Moore algebra
$\alpha\colon \Rdn(X)\rightarrow X$ is an affine function. For each
continuous map $f\colon X \rightarrow Y$, the function $\Rdn(f) \colon
\Rdn(X) \rightarrow \Rdn(Y)$ is affine.
\end{lem}

\proof
This follows from naturality of $\tau \colon \Dst U \Rightarrow
U\Rdn$. \qed

\auxproof{
To be precise, we have the following commuting diagrams.
$$\xymatrix@R-1pc{
\Dst U\Rdn(X)\ar[d]_{\tau}\ar[rr]^-{\Dst(U \alpha)} & &
   \Dst(UX)\ar[d]^{\tau} 
& &
\Dst U\Rdn(X)\ar[d]_{\tau}\ar[rr]^-{\Dst(U\Rdn(f))} & &
   \Dst U\Rdn(U)\ar[d]^{\tau} \\
U\Rdn^{2}(X)\ar[d]_{U(\mu)} & & U\Rdn(X)\ar[d]^{U(\alpha)} 
& &
U\Rdn^{2}(X)\ar[d]_{U(\mu)} & & U\Rdn^{2}(Y)\ar[d]^{U(\mu)} \\
U\Rdn(X)\ar[rr]_-{U(\alpha)} & & U(X)
& &
U\Rdn(X)\ar[rr]_-{U\Rdn(f)} & & U\Rdn(Y)
}$$
}



\begin{prop}
\label{EMHomProp}
Let $\alpha\colon\Rdn(X)\rightarrow X$ and $\beta\colon\Rdn(Y) \rightarrow Y$
be two Eilenberg-Moore algebras of the Radon monad $\Rdn$. A function
$f\colon X \rightarrow Y$ is an algebra homomorphism if and only if $f$ is
both continuous and affine.

As a result, the functor $\EM(\Rdn)\rightarrow \EM(\Dst) = \Conv$ from
Lemma~\ref{MonadFunctorLemma} is faithful, and an $\EM(\Dst)$ map comes
from an $\EM(\Rdn)$ map if and only if it is continuous.
\end{prop}

We shall follow the convention of writing $\Aff(X,Y)$ for the
homset of continuous and affine functions $X\rightarrow Y$.

\proof
Clearly, each algebra map is both continuous and affine. For the
converse, if $f\colon X\rightarrow Y$ is continuous, it is a map in
the category $\CH$ of compact Hausdorff spaces. Since it is affine,
both triangles commute in:
$$\xymatrix@R+0pc{
\Dst(UX)\ar@{ >->}[rr]^-{\tau}_-{\text{dense}}\ar[drr] & &
   \Rdn(X)\ar@<+1ex>[d]^{\beta\after\Rdn(f)}\ar@<-1ex>[d]_{f\after \alpha} \\
& & Y
}$$

\noindent Since $Y$ is Hausdorff, there is at most one such map. Therefore $f$ is an algebra map.
\qed

The category $\EM(\Rdn)$ of Eilenberg-Moore algebras of the Radon
monad may thus be understood as a suitable category of convex compact
Hausdorff spaces, with affine continuous maps between them. In the
next section, we see how to use a result from \cite{Swirszcz74} to
relate this to $\CCL$, which is a category of ``concrete'' convex
sets. Using this theorem, it will be shown that ``observability''
conditions like in~\cite[top of p. 169]{JacobsM12b} always hold for
algebras of $\Rdn$.

\subsection{\'Swirszcz's Theorem and Noncommutative \texorpdfstring{$C^*$}{C*}-algebras}
In this section we show that the Radon monad arises from an adjunction
in \cite{Swirszcz74} enabling us to use \'Swirszcz's theorem 3 from
that paper to show that the categories $\CCL$ and $\EM(\Rdn)$ are
equivalent, which we can then apply to represent noncommutative
$C^*$-algebras. The adjunction in question has $U \colon \CCL
\rightarrow \CH$ as the right adjoint, and the details of the
construction of the left adjoint are not given. In order to prove that
$\Rdn$ is the monad arising from this adjunction, we need to know its
unit and counit, so our next task is to define the left adjoint
explicitly. Of course, any other left adjoint will be naturally
isomorphic.

We begin as follows. We define $\Sw \colon \CH \rightarrow \CCL$ as
$\Sw = \Stat \after C$. Hence $\Rdn = U \after \Sw$.  To show that $\Sw$
is the left adjoint to $U$, we use the unit and counit definition of
an adjunction. We already know the unit, $\eta_X \colon X \rightarrow
U(\Sw(X))$, as we gave it when defining the unit of $\Rdn$. To define
the counit we use the notion of \emph{barycentre}.

We can understand the intuitive notion of barycentre by thinking of a Radon probability measure $\mu$ on the unit square $[0,1]^2$. If we wanted to find the centre of mass of $\mu$, which we shall call $b \in [0,1]^2$, we would take
$$\begin{array}{rclcrcl}
b_x 
& = &
\int\limits_{[0,1]^2} x \mathrm{d}\mu
& \qquad\mbox{and}\qquad &
b_y 
& = &
\int\limits_{[0,1]^2} y \mathrm{d}\mu
\end{array}$$

\noindent for the $x$ and $y$ coordinates. We can see that $x$ and $y$
are continuous affine functions from $[0,1]^2 \rightarrow \R$,
assigning each point to its $x$ and $y$ coordinate
respectively. Therefore we can rewrite the above as
$$\begin{array}{rclcrcl}
\int\limits_{[0,1]^2} x \mathrm{d}\mu 
& = &
x(b)
& \qquad\mbox{and}\qquad &
\int\limits_{[0,1]^2} y \mathrm{d}\mu 
& = &
y(b).
\end{array}$$

\noindent In monadic terms, this means that both projections $\pi_{1},
\pi_{2} \colon [0,1]^{2} \rightarrow [0,1]$ are maps of Eilenberg-Moore
algebras for the Radon monad, in the sense that the following diagram
commutes.
$$\xymatrix@R-.5pc{
\Rdn([0,1]^{2})\ar[d]_{\beta}\ar[rr]^-{\Rdn(\pi_{i})} & & 
   \Rdn([0,1])\ar[d]^{\alpha} \\
[0,1]^{2}\ar[rr]_-{\pi_i} & & [0,1]
}$$

\noindent We write $\alpha$ for the algebra $\nu \mapsto
\int\idmap\mathrm{d}\nu$, see also~\cite{Jacobs13a}, and $\beta$ for
the product algebra structure, given by $\mu \mapsto \tuple{\int
  \pi_{1}\mathrm{d}\mu, \int \pi_{2}\mathrm{d}\mu} = \tuple{\int
  x\mathrm{d}\mu, \int y\mathrm{d}\mu}$.

\auxproof{
The product algebra $\beta$ is obtained as:
$$\begin{array}{rcl}
\beta(\mu)
& = &
\tuple{\alpha(\Giry(\pi_{1})(\mu)), \alpha(\Giry(\pi_{2})(\mu))} \\
& = &
\tuple{\int \idmap \mathrm{d}\Giry(\pi_{1})(\mu),
   \int \idmap \mathrm{d}\Giry(\pi_{2})(\mu)} \\
& = &
\tuple{\int \pi_{1} \mathrm{d}\mu, \int \pi_{2} \mathrm{d}\mu}
   \qquad \mbox{by~\cite{Jacobs13a}}.
\end{array}$$
}

If we generalize $\pi_1$ and $\pi_2$ to arbitary real-valued continuous affine functions on $X$, and reinterpret Radon measures as functionals (as in the start of \S 5), we get the idea behind the following standard definition.

\begin{defi}
\label{BarycentreDefn}
If $X \in \CCL$ and $\phi \in \Sw(U(X))$, then a point $x \in X$ is a
\emph{barycentre} for $\phi$ if for all continuous affine functions
$f$ from $X \rightarrow \R$ we have that $\phi(f) = f(x)$. \qed
\end{defi}

The theorem that every $\phi$ has a barycentre when $X$ is a compact
subset of a locally convex space is standard and is proven in
\cite[proposition I.2.1 and I.2.2]{alfsen71}.

We will require the following important lemma, one of sevaral variants
of the Hahn-Banach separation lemma, and some of its corollaries,
which give an affine analogue of Urysohn's lemma for objects in
$\CCL$.

\begin{lem}
\label{ConvexSeparationLemma}
If $V$ is a locally convex topological vector space, $X$ a closed
convex subset and $Y$ a compact convex subset that is disjoint from
$X$, then there exists a continuous linear functional $\phi \colon V
\rightarrow \R$ and $\alpha \in \R$ such that $\phi(X) \subseteq
(\alpha, \infty)$ and $\phi(Y) \subseteq (-\infty, \alpha)$. \qed
\end{lem}
For proof, see either \cite[theorem IV.3.9]{Conway90} or \cite[II.4.2 corollary 1]{Schaefer66}.

\begin{cor}
\label{ConvexSepCorr}
Let $(K,V) \in \Obj(\CCL)$. In the following $X,Y$ will be arbitrary closed disjoint convex subsets of $K$, $x,y$ arbitrary distinct points of $K$. 
\begin{enumerate}[label=(\roman*)]
\item There is a $\phi \in \Aff(K, \R)$ and an $\alpha \in \R$ such that $\phi(X) \subseteq (\alpha, \infty)$ and $\phi(Y) \subseteq (-\infty, \alpha)$. 
\item There is a $\phi \in \Aff(K, \R)$ such that $\phi(x) \neq \phi(y)$.
\item There is a $\phi \in \CCL(K,[0,1])$ and an $\alpha \in \R$ such that $\phi(X) \subseteq (\alpha, 1]$ and $\phi(Y) \subseteq [0, \alpha)$.
\item There is a $\phi \in \CCL(K,[0,1])$ such that $\phi(x) \neq \phi(y)$.
\end{enumerate}
\end{cor}

\proof\hfill
\begin{enumerate}[label=(\roman*)]
\item Apply Lemma \ref{ConvexSeparationLemma} to obtain $\phi' \colon
  V \rightarrow \R$ separating $X$ from $Y$. Since $K$ has the
  subspace topology, $\phi = \phi'|_K$ is continuous, and since
  $\phi'$ is linear, $\phi$ is affine, hence $\phi \in \Aff(K,\R)$. We
  also keep the properties that $\phi(X) \subseteq (\alpha, \infty)$
  and $\phi(Y) \subseteq (-\infty, \alpha)$.
\item This follows directly from (i), using the fact that points are compact and convex.
\item We use (i) and obtain $\phi' \in \Aff(K,\R)$ and $\alpha' \in \R$. Since the image of a compact space is compact, and a compact subset of $\R$ is closed and bounded, the numbers
\begin{align*}
\beta_{\uparrow} &= \sup \phi'(K) &
\beta_{\downarrow} &= \inf \phi'(K)
\end{align*}
exist, and $\phi'$ can be considered as an affine continuous map $K \rightarrow [\beta_{\downarrow}, \beta_{\uparrow}]$. We define
\[
\phi(k) = \frac{\phi(k)-\beta_{\downarrow}}{\beta_{\uparrow}-\beta_{\downarrow}}
\]
if $\beta_{\uparrow} \neq \beta_{\downarrow}$, otherwise we define it without dividing by anything, though this can only happen if one of $X$ or $Y$ is empty. The image of $\phi$ is contained in $[0,1]$, and $\phi$ is affine and continuous, being the composition of affine and continuous maps. We define
\[
\alpha = \frac{\alpha' - \beta_{\downarrow}}{\beta_{\uparrow} - \beta_{\downarrow}}
\]
again not doing the division if it is zero. We have that $\phi(X) \subseteq (\alpha, \infty)$, and since the image of $\phi$ is contained in $[0,1]$, this implies $\phi(X) \subseteq (\alpha, 1]$. The proof that $\phi(Y) \subseteq [0, \alpha)$ is similar.
\item This is proven using (iii), again using the fact that points are closed, convex sets. \qed
\end{enumerate}

\noindent Using the properties proven above, we can start to define the counit of the adjunction.

\begin{lem}
\label{SwirszczCounitWellAffLemma}\hfill
\begin{enumerate}[label=(\roman*)]
\item For every $\phi \in \Sw(U(X))$ the barycentre is unique. The
  function $\varepsilon_X \colon \Sw(U(X)) \rightarrow X$ mapping
  $\phi$ to its barycentre is well defined.
\item This $\varepsilon_X$ is an affine map.
\end{enumerate}
\end{lem}
\proof\hfill
\begin{enumerate}[label=(\roman*)]
\item We show the barycentre is unique as follows. Let $(V, X)$ be an object of $\CCL$, $V$ being the locally convex space and $X$ the compact convex subset. Let $x, x' \in X$ be barycentres of $\phi \in \Sw(U)$. Suppose for a contradiction that $x \neq x'$. By corollary \ref{ConvexSepCorr} (ii), there is an $f \in \Aff(X,\R)$ such that $f(x) \neq f(x')$. Since $x$ and $x'$ are both barycentres of $\phi$,
\[
f(x) = \phi(f) = f(x')
\]
a contradiction. So we have $x = x'$. Therefore $\varepsilon_X$ is well-defined, at least as a function between sets. 
\item To show that $\varepsilon_X$ is affine, consider two Radon
  measures $\phi, \psi \in \Sw(U(X))$, such that $\varepsilon_X(\phi)
  = x$ and $\varepsilon_X(\psi) = y$, i.e. these are the
  barycentres. To show that $\varepsilon_X(\alpha \phi + (1-\alpha)
  \psi) = \alpha \varepsilon_X(\phi) + (1-\alpha)
  \varepsilon_X(\psi)$, we will show that $\alpha x + (1 - \alpha) y$
  is the barycentre of $\alpha \phi + (1- \alpha) \psi$. Given an
  continuous affine function $f \colon X \rightarrow \R$, we have
\[
(\alpha \phi + (1 - \alpha) \psi)(f) = \alpha \phi(f) + (1 - \alpha) \psi(f) = \alpha x + (1 - \alpha) y
\]
so $\varepsilon_X$ is affine. \qed
\end{enumerate}

\begin{lem}
\label{SwirszczCounitContinuousLemma}
The barycentre map $\varepsilon_X$ is continuous, hence a map in $\CCL$.
\end{lem}
\proof
We now show that $\varepsilon_X$ is continuous. We use the filter-theoretic definition of continuity. Given $\phi \in \Sw(U(X))$, with barycentre $x$, we want to show that for every neighbourhood $V$ of $x$, there is a neighbourhood $U$ of $\phi$ such that $\varepsilon_X(U) \subseteq V$. It suffices to prove this for a chosen set of basic neighbourhoods, so we choose open neighbourhoods for $X$ and for $\Sw(U(X))$ we choose finite intersections of elements of the following subbasis of closed neighbourhoods:
\[
U_{f,\alpha,\epsilon} = \{ \psi \in \Sw(U(X)) \mid | \psi(f) - \alpha | \leq \epsilon \}
\]
where $f \in C(U(X))$, $\alpha \in \R$ and $\epsilon \in (0,\infty)$.

We find the neighbourhood of $\phi$ using a compactness argument. 

Consider the following subset of $X$.
\[
\bigcap_{\substack{f \in \Aff(X,\R) \\ \epsilon > 0}} \overline{\varepsilon_X(U_{f,f(x),\epsilon})}
\]
Since $\phi \in U_{f,f(x),\epsilon}$ for all values of $f$ and $\epsilon$, we have that $x$ is in this intersection. We will show that
\begin{equation}
\label{BigIntersectionIsPoint}
\bigcap_{\substack{f \in \Aff(X,\R) \\ \epsilon > 0}} \overline{\varepsilon_X(U_{f,f(x),\epsilon})} = \{ x \}
\end{equation}
As we already know $x$ is an element of the left hand side, we will show that if $x' \in X$ and $x' \neq x$, then $x'$ is not an element of the left hand side. So since $x \neq x'$, by Corollary \ref{ConvexSepCorr}(ii) there is an $f \in \Aff(X, \R)$ such that $f(x) \neq f(x')$. We let 
\begin{equation}
\label{EpsilonDefInequality}
\epsilon = \frac{|f(x) - f(x')|}{3} > 0
\end{equation}
We show that $x' \not\in \overline{\varepsilon_X(U_{f,f(x),\epsilon})}$ and therefore is not in \eqref{BigIntersectionIsPoint} by showing there is an open set containing $x'$ that is disjoint from $\varepsilon_X(U_{f,f(x),\epsilon})$. The open set we choose is
\[
f^{-1}((f(x')-\epsilon, f(x') + \epsilon))
\]
which is open because $f$ is continuous. Assume for a contradiction that there is some $x'' \in f^{-1}((f(x')-\epsilon, f(x') + \epsilon)) \cap \varepsilon_X(U_{f,f(x),\epsilon})$. This means that 
\begin{equation}
\label{FXsEpsilonInequality}
|f(x') - f(x'')| < \epsilon
\end{equation} and there is some $\psi \in U_{f,f(x),\epsilon}$ of which $x''$ is the barycentre, \emph{i.e.} for all $g \in \Aff(X, \R)$ $\psi(g) = g(x'')$. Therefore it must be the case that $\psi(f) = f(x'')$, and so the inequality deriving from $\psi \in U_{f,f(x),\epsilon}$, which is $|\psi(f) - f(x) | \leq \epsilon$ becomes $|f(x'') - f(x)| \leq \epsilon$. If we combine this with \eqref{FXsEpsilonInequality} and use the triangle inequality, we get $|f(x') - f(x)| \leq 2\epsilon$, which contradicts $|f(x)-f(x')| \geq 3\epsilon$ from \eqref{EpsilonDefInequality}. Therefore the assumption that $x''$ could exist is wrong, so $x'$ is in an open set outside $\varepsilon_X(U_{f,f(x),\epsilon})$, and hence $x' \not\in \overline{\varepsilon_X(U_{f,f(x),\epsilon})}$. This establishes that \eqref{BigIntersectionIsPoint} is the case. 

Now consider $X \setminus V$, which is a closed set that does not contain $x$, since $V$ is an open neighbourhood of $x$. We therefore have
\[\emptyset = (X \setminus V) \cap \bigcap_{\substack{f \in \Aff(X,\R) \\ \epsilon > 0}} \overline{\varepsilon_X(U_{f,f(x),\epsilon})}\;
 = \bigcap_{\substack{f \in \Aff(X,\R) \\ \epsilon > 0}} (X \setminus V) \cap \overline{\varepsilon_X(U_{f,f(x),\epsilon})}
\]

The right hand side is a family of closed subsets of a compact space with empty intersection. Therefore there is a finite subfamily also having empty intersection. We use the numbers $i \in \{1, \ldots, n\}$ as an index set, and take $\{\epsilon_i\}$, $\{f_i\}$ such that we have
\[\emptyset = \bigcap_{i = 1}^n (X \setminus V) \cap \overline{\varepsilon_X(U_{f_i,f_i(x),\epsilon_i})}\;
 = (X \setminus V) \cap \bigcap_{i=1}^n \overline{\varepsilon_X(U_{f_i,f_i(x),\epsilon_i})}
\]
Therefore we have
\[
\varepsilon_X\left(\bigcap_{i=1}^n U_{f_i,f_i(x),\epsilon_i}\right) \subseteq \bigcap_{i=1}^n \varepsilon_X(U_{f_i,f_i(x),\epsilon_i}) \subseteq \bigcap_{i=1}^n \overline{\varepsilon_X(U_{f_i,f_i(x),\epsilon_i})} \subseteq V
\]
Since $V$ was an arbitrary open neighbourhood of $\varepsilon_X(\phi)$, we have that $\varepsilon_X$ is continuous at $\phi$. Since the choice of $\phi$ was arbitrary, $\varepsilon_X$ is continuous. 
\qed

\begin{lem}
\label{SwirszczCounitNaturalLemma}
The family $\{\varepsilon_X\}$ defines a natural transformation $\varepsilon: \Sw \after U \Rightarrow \Idmap$.
\end{lem}
\proof
We must show that
\[
\xymatrix@C+1pc{
\Sw(U(X)) \ar[r]^-{\epsilon_X} \ar[d]_{\Sw(U(f))} & X \ar[d]^f \\
\Sw(U(Y)) \ar[r]_-{\epsilon_Y} & Y
}
\]
Suppose that $\phi \in \Sw(U(X))$ and $\varepsilon_X(\phi) = x$, i.e. $x$ is the barycentre of $\phi$. It suffices to show that $f(x)$ is the barycentre of $\Sw(U(f)(\phi)$. Let $h \in C(Y)$, and we have by definition that
\[
\Sw(U(f))(\phi)(h) = \phi(h \after f)
\]

We want to show that if $h$ is affine, then $\Sw(U(f))(\phi)(h) = h(f(x))$, as this would show $f(x)$ is the barycentre. Since $h \after f$ is the composite of continuous, affine functions, it is also continuous and affine, and so, using that $x$ is the barycentre of $\phi$, we have that $\phi(h \after f) = (h \after f)(x) = h(f(x))$, which is what we were required to prove. 
\qed

Taken together, the preceding three lemmas define the counit. We can now move on to showing that this is actually an adjunction.

\begin{thm}
\label{SwirszczAdjunctionTheorem}
The functor $\Sw \colon \CH \rightarrow \CCL$ is the left adjoint to
$U \colon \CCL \rightarrow \CH$
\end{thm}
\proof
We show that the unit-counit diagrams commute.

First we must show that the following commutes:
\[
\xymatrix{
U Y \ar[r]^-{\eta_{U Y}} \ar[rd]_{\id_{U Y}} & U(\Sw(U(Y))) \ar[d]^{U \varepsilon_Y} \\
 & U Y
}
\]
In other words, we must show that for all $y \in U Y$, $y$ is the
barycentre of $\eta_{U Y}(y)$. Using the definition of $\eta$, we have
that for any affine continuous function $f \colon X \rightarrow \R$
that
\[
\eta_{U Y}(x)(f) = f(x)
\]
because that is already true for all continuous functions $f \in C(X)$. Therefore $x$ is the barycentre of $\eta_{U Y}(x)$, and so the diagram commutes.

The second diagram we must consider is the following:
\[
\xymatrix{
\Sw(X) \ar[r]^-{\Sw(\eta_X)} \ar[rd]_{\id_{\Sw(X)}} & \Sw(U(\Sw(X))) \ar[d]^{\epsilon_{\Sw(X)}} \\
 & \Sw(X)
}
\]
This time, we need to show that $\phi \in \Sw(X)$ is the barycentre of
the measure $\Sw(\eta_X)(\phi)$. So consider an affine continuous
function $k \colon \Sw(X) \rightarrow \R$. We want to show that
$\Sw(\eta_X)(\phi)(k) = k(\phi)$ for all $\phi \in \Sw(X)$. To do
this, we use Lemma \ref{DensityLemma}. We show the diagram commutes on
the convex combinations of extreme points, and since this is a dense
subset, the diagram commutes by continuity. So let $\{x_1, \ldots
x_n\}$ be a finite subset of $X$, and
\[
\sum_{i=1}^n \alpha_i\eta_X(x_i)
\]
a finite convex combination of extreme points of $\Sw(X)$. Now
\begin{align*}
\Sw(\eta_X)\left(\sum_{i=1}^n \alpha_i\eta_X(x_i)\right)(k) &= \left(\sum_{i=1}^n\alpha_i \eta_X(x_i) \right)(k \after \eta_X) \\
 &= \sum_{i=1}^n \alpha_i \eta_X(x_i)(k \after \eta_X) \\
 &= \sum_{i=1}^n \alpha_i k(\eta_X(x_i)) \\
 &= k\left(\sum_{i=1}^n(\eta_X(x_i))\right)
\end{align*}
with the last step holding because $k$ is an affine function. 

As explained before, this shows $\Sw(\eta_X)(\phi)(k) = k(\phi)$ for all $\phi \in \Sw(X)$, and hence the diagram commutes.
Thus we have that $\Sw$ is the left adjoint to $U$. 
\qed

Now that we have defined the adjunction $\Sw \dashv U$, we can move on
to proving that $\Rdn$ is not only the same functor as the monad
derived from $\Sw \dashv U$ but also the same as a monad. In order to
do this, we require a few lemmas concerning the definition of $\mu$ we
gave at the start of Section \ref{ContProbSec}. The map $\mu$ was
defined using $\lam{h}{h(v)}$. Since we need to prove certain
properties about it, we give this map a name, and generalize it
somewhat for later use. If $A$ is a (possibly noncommutative)
$C^*$-algebra, we define
$$\xymatrix{
A\SA\ar[r]^-{\zeta_A} & \Aff(\Stat(A), \R) 
\qquad\mbox{as}\qquad
\zeta_A(a)(\phi) = \phi(a).
}$$

\noindent In the special case we had earlier, we were using
$\zeta_{C(X)}$ for a compact Hausdorff space $X$, since $C(X)\SA =
C_\R(X)$, the real-valued functions. We can see that
\begin{equation}
\label{MuDefTwo}
\mu_X(g)(v) = g(\zeta_{C(X)}(v)).
\end{equation}

\begin{lem}
\label{ZetaIsoLemma}
The map $\zeta_A$ is a bijection between $A\SA$ and $\Aff(\Stat(A),\R)$. $\zeta_{C(X)}$ is a bijection between $C_\R(X)$ and $\Aff(\Sw(X), \R)$. In fact, the bijection is an isomorphism of ordered $\R$-vector spaces with unit, taking these to be defined pointwise on $\Aff(\Stat(A), \R)$. \qed
\end{lem}
The proof can be found in \cite[Proposition 2.3]{Alfsen01}. It was
originally proved by Kadison \cite[Lemma 4.3, Remark 4.4]{Kadison65}
and is often stated for complete order-unit spaces (such as in
\cite[Theorem II.1.8]{alfsen71}), though it was originally intended
for use with $C^*$-algebras, as here.

\auxproof{
\begin{itemize}
\item $\zeta_X(f)$ is affine:
\begin{align*}
\zeta_X(f)(\alpha \phi + \beta \psi) &= (\alpha \phi + \beta \psi)(f) \\
 &= \alpha \phi(f) + \beta \psi(f) \\
 &= \alpha \zeta_X(f)(\phi) + \beta \zeta_X(f)(\psi)
\end{align*}
\item $\zeta_X(f)$ is continuous:

Let $\phi_\alpha$ be a net converging to $\phi$ in the weak-* topology
on $\Sw(X)$, which is to say that for all $f \in C(X)$,
$\phi_\alpha(f)$ converges to $\phi(f)$. Then for all $f \in C(X)$, we
have that
\begin{align*}
\lim_\alpha \zeta_X(f)(\phi_\alpha) &= \lim_\alpha \phi_\alpha(f) \\
 &= \phi(f) \\
 &= \zeta_X(f)(\phi)
\end{align*}
hence $\zeta_X(f)$ is continuous as it preserves limits of nets.

\item $\zeta_X$ is injective:

Suppose $\zeta(f) = \zeta(g)$, for $f,g \in C(X)$. Then for all $\phi \in \Sw(X)$, $\phi(f) = \phi(g)$. Therefore we have that for all $x \in X$
\[
\eta_X(x)(f) = \eta_X(x)(g)
\]
hence $f(x) = g(x)$, and so $f = g$. 

\item $\zeta_X$ is surjective:

Let $f \in \Aff(\Sw(X), \R)$, which is therefore a continuous function
$f \colon U(\Sw(X)) \rightarrow \R$. We may take $f \after \eta_X$,
which is a continuous function $X \rightarrow \R$, i.e. an element of
$C_\R(X)$. So we must show that $\zeta_X(f \after \eta_X) = f$. We do
this by showing $\zeta(f \after \eta_X)(\phi) = f(\phi)$ where $\phi$
is an element of the dense subset of finite convex combinations of
$\eta_X(x)$ measures. Let $\phi = \sum_{i=1}^n \alpha_i \eta(x_i)$ be
such a convex combination.
\begin{align*}
\zeta(f \after \eta_X)(\phi) &= \zeta(f \after \eta_X)\left(\sum_{i=1}^n \alpha_i \eta(x_i) \right) \\
 &= \left(\sum_{i=1}^n \alpha_i \eta(x_i) \right)(f \after \eta_X) \\
 &= \sum_{i=1}^n \alpha_i \eta(x_i)(f \after \eta_X) \\
 &= \sum_{i=1}^n \alpha_i f(\eta_X(x_i)) \\
 &= f\left(\sum_{i=1}^n \alpha_i \eta_X(x_i) \right) \\
 &= f(\phi)
\end{align*}
hence this last statement follows by lemma \ref{DensityLemma}.
\end{itemize}
}

\begin{thm}
\label{SwirszczRadonMonadTheorem}
The monad $:\CH \rightarrow \CH$ given by $\Sw \dashv U$ is the Radon monad $\Rdn$.
\end{thm}
\proof
We have by definition that $\Rdn = U \Sw$ and $\eta = \eta$. Therefore we only need to show that $\mu = U \varepsilon \Sw$. What we need to show then, is that if $X$ is a compact Hausdorff space and $\phi \in \Sw(U(\Sw(X)))$, then $\mu(\phi)$ is the barycentre of $\phi$. That is to say, for all $f \in \Aff(\Sw(X), \R)$, $\phi(f) = f(\mu_X(\phi))$. Using Lemma \ref{ZetaIsoLemma}, we reduce to showing that for all $f \in C_\R(X)$, we have $\phi(\zeta_X(f)) =  \zeta_X(f)(\mu_X(\phi))$. Using \eqref{MuDefTwo}, we have
\[
\zeta_X(f)(\mu_X(\phi)) = \mu_X(\phi)(f) = \phi(\zeta_X(f))
\]
as required. \qed

\auxproof{
Here are some of the definitions we will need.

The map $\gamma_X \colon \E(X) \rightarrow \Sw(\beta(X))$ is given by
\[
\gamma_X(\phi \colon [0,1]^X \rightarrow [0,1])(f \colon \beta(X)
\rightarrow \C) = \phi(\theta_{\ell^\infty(X)}(f \after \eta^\beta_X))
\]
where $\theta_{\ell^\infty(X)}$ is the map extending effect algebra maps $[0,1]^X \rightarrow [0,1]$ to states $\ell^\infty(X) \rightarrow \C$, and $\eta^\beta$ is the unit of the Stone-\v{C}ech adjunction, mapping elements of $X$ to principal ultrafilters. We should show this is a monad morphism, and since it is an isomorphism, it is a monad isomorphism.  

The $\eta$ and $\mu$ for $\E$ are given as follows, from \cite[page
  157]{JacobsM12b}. The letters are $p \in [0,1]^X$, $h \colon
[0,1]^{\E(X)} \rightarrow [0,1]$ and $k \in \E(X)$.
\begin{align*}
\eta_X(x)(p) &= p(x) \\
\mu_X(h)(p) &= h(\lam{k}{k(p)})
\end{align*}
we have something similar to $\zeta_X$ above. We could call it $\xi_X$ and prove that it gives all affine functions on $\Sw(\beta(X))$, starting with bounded functions on $X$. 

\cite[page 150]{JacobsM12b} gives the unit and counit of a composite adjunction. These are given in this case by
\begin{align*}
\eta^{\Sw\beta}_X &= V \eta^\Sw_{\beta(X)} \after \eta^\beta_X \\
\varepsilon^{\Sw\beta}_Y &= \varepsilon^\Sw_Y \after \Sw \varepsilon_{UY} \\
\mu_X &= U_\Sets \varepsilon_{\Sw \beta X} = U_\Sets (\varepsilon^\Sw_{\Sw \beta X} \after \Sw(\varepsilon^\beta_{U\Sw\beta X}))
\end{align*}
}

\begin{thm}[\'Swirszcz's theorem]
\label{SwirszczTheorem}
The forgetful functor $U \colon \CCL \rightarrow \CH$ is monadic,
\emph{i.e.} $\CCL \simeq \EM(U \after \Sw)$. By Theorem
\ref{SwirszczRadonMonadTheorem}, $\CCL \simeq \EM(\Rdn)$. \qed
\end{thm}
This comes from \cite[Theorem 3]{Swirszcz74}. A proof not using any monadicity theorems can be found in \cite[Proposition 7.3]{Semadeni73}.

\subsubsection{Non-commutative $C^*$-algebras and $\EM(\Rdn)$}
In the following section we shall show that the category $\CstarPU$
embeds fully and faithfully in $\EM(\Rdn)$. To do this, we use the
fact that $\EM(\Rdn) \simeq \CCL$, and also the functor $\Stat \colon
\CstarPU \rightarrow \CCL$.

We begin with a standard separation result from the theory of $C^*$-algebras.

\begin{lem}
\label{StateSepLemma}
If $A$ is a $C^*$-algebra, and $a,b \in A$, then
\[
\phi(a) = \phi(b)
\]
for all $\phi \in \Stat(A)$ implies $a = b$. In other words, $A$ is separated by its states, or $A$ has ``sufficiently many states''.
\end{lem}
\proof
In \cite[theorem 4.3.4 (i)]{KadisonR83} we have that if $\phi(a) = 0$ for all $\phi \in \Stat(A)$, then $a = 0$. We simply apply this to $a-b$. 
\qed

On the set $\Aff(X, \C)$, for $X \in \Obj(\CCL)$, we can define a $\C$-vector space structure, a positive cone, and a distinguished unit, simply by using the fact that $\C$ has these things and defining them pointwise. The positive cone is $[0, \infty) \subseteq \C$ and the unit is $1$. Given these definitions, we can prove the complexification of Lemma \ref{ZetaIsoLemma}.

\begin{lem}
\label{ComplexKadisonLemma}
For each $C^*$-algebra $A$, the map $\xi_A \colon A \rightarrow
\Aff(\Stat(A),\C)$, defined as
\[
\xi_A(a)(\phi) = \phi(a)
\]
is an isomorphism of complex vector spaces preserving the positive
cone and unit in both directions.
\end{lem}
\proof
First we show that the map $\xi_A$ is $\C$-linear and preserves $^*$. For $\C$-linearity, let $z \in \C$, $\phi \in \Stat(A)$ and $a \in A$. Then
$$\begin{array}{rcccccl}
\xi_A(za)(\phi)
& = &
\phi(za)
& = &
z\phi(a) 
& = &
z \xi_A(a)(\phi),
\end{array}$$

\noindent so $\xi_A(za) = z\xi_A(a)$. 

To show that it preserves $^*$, where for $f \in \Aff(\Stat(A),\C)$, $f^*$ is calculated pointwise, we use the fact that every positive linear functional on $A$, and hence every state, is self-adjoint, as described in Lemma \ref{CstarMapLem}, \emph{i.e.} $\phi(a^*) = \overline{\phi(a)}$. 

Thus we have
$$\begin{array}{rcccccccl}
\xi_A(a^*)(\phi)
& = &
\phi(a^*)
& = &
\overline{\phi(a)}
& = &
\overline{\xi_A(a)(\phi)}
& = &
\xi_A(a)^*(\phi).
\end{array}$$

\noindent and so $\xi_A(a^*) = \xi_A(a)^*$. 

From Lemma~\ref{ZetaIsoLemma} we have that $\xi$ restricts to an
isomorphism $\zeta \colon A\SA \cong \Aff(\Stat(A), \R)$ as an ordered
vector space with unit. We extend this to complex numbers as
follows. Given $a \in A$, we can define its real and imaginary parts
as
\begin{align*}
\Re(a) &= \frac{a + a^*}{2} &
\Im(a) &= \frac{a - a^*}{2i}
\end{align*}
and we see that $\Re(a) + i \Im(a) = a$. Similarly, using pointwise complex conjugation as $^*$, we can define real and imaginary parts of an affine continuous map from $\Stat(A) \rightarrow \C$, and the self-adjoint elements are maps $\Stat(A) \rightarrow \C$. Since we know that $\eta_X$ has an inverse for self-adjoint elements, we can define the inverse as
\[
\xi_A^{-1}(f + ig) = \xi_A^{-1}(f) + i\xi_A^{-1}(g)
\]
where $f,g$ are self-adjoint.

We show this is the inverse of $\xi_A$. For one way
\begin{align*}
\xi_A(\xi_A^{-1}(f + ig)) &= \xi_A(\xi_A^{-1}(f) + i\xi_A^{-1}(g)) \\
 &= \xi_A(\xi_A^{-1}(f)) + i\xi_A(\xi_A^{-1}(g)) \\
 &= f + ig .
\end{align*}
For the other way, with $a,b \in A^{\mathrm{sa}}$,
\begin{align*}
\xi_A^{-1}(\xi_A(a + ib)) &= \xi_A^{-1}(\xi_A(a) + i\xi_A(b)) \\
 &= \xi_A^{-1}(\xi_A(a)) + i\xi_A^{-1}(\xi_A(b)) \\
 &= a + ib ,
\end{align*}
where the definition of $\xi_A^{-1}$ can be applied since $\xi_A$ preserves $^*$ and hence preserves self-adjointness, so $\xi_A(a)$ and $\xi_A(b)$ are both self-adjoint.
\qed

We will require the following fact in a moment.
\begin{lem}
\label{StateSpaceLemma4}
If $B$ is a $C^*$-algebra, $b' \in \Aff(\Stat(B), \C)$, then for all $\phi \in \Stat(B)$
\[
\phi(\xi_B^{-1}(b')) = b'(\phi).
\]
\end{lem}
\proof
By Lemma \ref{ComplexKadisonLemma}, we have that there is some $b \in B$ such that $b' = \xi_B(b)$. Then we have
\[
\phi(\xi_B^{-1}(b'))
 = 
\phi(\xi_B^{-1}(\xi_B(b)))
 = 
\phi(b)
 = 
\xi_B(b)(\phi)
 = 
b'(\phi).\eqno{\qEd}
\]\smallskip

\noindent We can now prove that $\Stat$ is full and faithful, and hence
$\op{(\CstarPU)}$ embeds fully in $\EM(\Rdn)$.

\begin{thm}
\label{StatFullFaithfulThm}
The state space functor $\Stat \colon \op{(\CstarPU)} \rightarrow
\CCL$ is full and faithful.
\end{thm}

\proof\hfill
\begin{itemize}
\item For faithfulness, suppose we have $f,g \colon A \rightarrow B$
  in $\CstarPU$, such that $\Stat(f) = \Stat(g)$. We have that
  $\Stat(f)(\phi) = \Stat(g)(\phi)$ for all $\phi \in \Stat(B)$,
  which, expanding the definitions, gives that $\phi \after f = \phi
  \after g$ for all $\phi \in \Stat(B)$. Now, we have that for all $a
  \in A$ and $\phi \in \Stat(B)$, that $\phi(f(a)) = \phi(g(a))$. By
  Lemma \ref{StateSepLemma}, we have that for all $a \in A$, $f(a) =
  g(a)$, and therefore $f = g$.

\item For fullness, let $g \colon \Stat(B) \rightarrow \Stat(A)$ be an
  affine, continuous map. We must find a map $f \colon A \rightarrow
  B$ such that $\Stat(f) = g$. We take the map $f = \xi_B^{-1} \after
  \Aff(g, \C) \after \xi_A \colon A \rightarrow B$. First we must
  prove this map is positive, $\C$-linear, and unital. We know from
  Lemma \ref{ComplexKadisonLemma} that, being isomorphisms, $\xi_A$
  and $\xi_B^{-1}$ are $\C$-linear (with the pointwise structure on
  $\Aff(\Stat(A), \C)$) and preserve the positive cone and
  unit. Therefore we only need to show that $\Aff(g, \C)$ has these
  properties to verify them for $f$. For $\C$-linearity, let $a_1, a_2
  \in \Aff(\Stat(A), \C)$, and $z_1, z_2 \in \C$. Then for each $\phi
  \in \Stat(B)$
\begin{align*}
\Aff(g,\C)(z_1a_1 + z_2a_2)(\phi) &= ((z_1a_1 + z_2a_2) \after g)(\phi) \\
 &= (z_1a_1+z_2a_2)(g(\phi)) \\
 &= z_1a_1(g(\phi)) + z_2a_2(g(\phi)) \\
 &= z_1\Aff(g,\C)(a_1)(\phi) + z_2\Aff(g,\C)(a_2)(\phi) \\
 &= (z_1\Aff(g,\C)(a_1) + z_2\Aff(g,\C)(a_2))(\phi),
\end{align*}
and so 
\[
\Aff(g,\C)(z_1a_1 + z_2a_2) = z_1\Aff(g,\C)(a_1) + z_2\Aff(g,\C)(a_2),
\]
which is to say, $\Aff(g,\C)$ is $\C$-linear.

The unit of $\Aff(\Stat(A),\C)$ is given by the function $1 \colon
\Stat(A) \rightarrow \C$ that maps every element of $\Stat(A)$ to $1
\in \C$. We must show that $\Aff(g,\C)$ preserves this unit. Given
$\phi \in \Stat(B)$, we have
\[
\Aff(g,\C)(1)(\phi) = 1(g(\phi)) = 1,
\]
so $\Aff(g, \C)(1)$ takes the value $1 \in \C$ for all $\phi \in
\Stat(B)$, and hence it is the unit in $\Aff(\Stat(B),\C)$. 

The positive elements of $\Aff(\Stat(A),\C)$ are given by functions whose image is contained in the positive reals, $[0,\infty) \subseteq \C$. We need to show that if $a \in \Aff(\Stat(A),[0,\infty))$, then so is $\Aff(g,\C)(a)$. This is easily accomplished as before. If $\phi \in \Stat(B)$, then
\[
\Aff(g,\C)(a)(\phi) = (a\after g)(\phi) = a(g(\phi)).
\]
Since $g(\phi) \in \Stat(A)$, we have that $a(g(\phi)) \in [0,\infty)$ by the assumption on $a$, and so $\Aff(g,\C)(a)$ is a positive element of $\Aff(\Stat(B),\C)$. All these conditions, taken together, show that $f$ is a $\CstarPU$ map from $A$ to $B$.

Now we show that $\Stat(f) = g$. Let $\phi \in \Stat(B)$ and $a \in A$. Then
\begin{align*}
\Stat(f)(\phi)(a) &= \Stat(\xi_B^{-1} \after \Aff(g, \C) \after \xi_A)(\phi)(a)\\
 &= (\phi \after \xi_B^{-1} \after \Aff(g, \C) \after \xi_A)(a) \\
 &= \phi(\xi_B^{-1}(\Aff(g, \C)(\xi_A(a)))) \\
 &= \phi(\xi_B^{-1}(\xi_A(a) \after g)),
\end{align*}
applying Lemma \ref{StateSpaceLemma4}, we continue
\begin{align*}
\Stat(f)(\phi)(a) &= (\xi_A(a) \after g)(\phi) \\
 &= \xi_A(a)(g(\phi)) \\
 &= g(\phi)(a).
\end{align*}
Since this holds for all $\phi$ and $a$, we have the required equality $\Stat(f) = g$, proving $\Stat$ is full. \qed
\end{itemize}

\noindent Alfsen, Hanche-Olsen and Shultz have characterized the essential image
of $\Stat$ \cite[Corollary 8.6]{Alfsen80}. We do not give the
characterization here as it involves many further definitions. Since
there are PU-maps that are not completely positive, $\Stat$ is not a
full functor when restricted to $\CstarCPU$. In fact, whether a map is
completely positive or not depends on the orientation (in the sense of
\cite{Alfsen80}) and cannot be defined purely from the $\EM(\Rdn)$
structure of the state space. This can be seen by the fact that the
transpose map, the archetypal positive but not completely positive
map, is self-inverse, and hence an isomorphism as a PU map, and so by
the above result defines an isomorphism in $\EM(\Rdn)$ on the state
space.



\section{States and effects}

We start with a simple observation.

\begin{lem}
\label{UnitRalgLem}
The unit interval $[0,1]$ is a compact convex subset of the locally convex space $\R$, and therefore carries a $\Rdn$-algebra structure by Theorem \ref{SwirszczTheorem}. The algebra map $\Rdn([0,1]) \rightarrow [0,1]$ maps each measure to its mean value.

For an arbitrary $\Rdn$-algebra $X$, the homset of algebra maps:
$$\begin{array}{rcl}
\EM(\Rdn)\big(X, [0,1]\big)
& = &
\Aff(X, [0,1])
\end{array}$$

\noindent is an effect module, with pointwise operations. Recall
from Proposition~\ref{EMHomProp} that this homset is the affine
and continuous functions $X \rightarrow [0,1]$. Taken all together, we have defined a functor $\Aff(-,[0,1]) \colon \EM(\Rdn) \rightarrow
\op{\EMod}$. \qed
\end{lem}


In~\cite{JacobsM12b} it is shown that for an effect module $M$, the
homset $\EMod(M,[0,1])$ is a convex compact Hausdorff space. In fact,
it carries an $\Rdn$-algebra structure:
$$\xymatrix@R-2pc{
\Rdn\big(\EMod(M,[0,1])\big)\ar[rr]^-{\alpha_M} & & \EMod(M,[0,1]) \\
h\ar@{|->}[rr] & & \lamin{x}{M}{h(\ev_{x})}
}$$

\noindent where $\ev_{x} = \lam{v}{v(x)} \colon C\big(\EMod(M,[0,1])\big)
\rightarrow \C$. For each map of effect modules $f\colon M \rightarrow
M'$ one obtains a map of $\Rdn$-algebras $(-) \after f \colon
\EMod(M',[0,1]) \rightarrow \EMod(M,[0,1])$. We thus obtain the
following situation:
\begin{equation}
\label{KlRSETDiag}
\vcenter{\xymatrix{
\op{\EMod}\ar@/^1.5ex/[rr]^-{\EMod(-,[0,1])} & \top & \EM(\Rdn)\ar@/^1.5ex/[ll]^-{\Aff(-,[0,1])} & \op{\EMod} \ar@/^1.5ex/[rr]^-{\EMod(-,[0,1])} & \top & \EM(\Rdn) \ar@/^1.5ex/[ll]^-{\Aff(-,[0,1])} \\
& \Kl(\Rdn)\ar[ul]^{\Cont(-,[0,1])\quad}\ar[ur] & & & \op{(\CstarPU)} \ar[ul]^{[0,1]_{(-)}} \ar[ur]_\Stat &
}}
\end{equation}

\noindent Such diagrams appear in~\cite{Jacobs12} as a categorical
representation of the duality between states and effects, with the
Schr\"odinger picture on the right vertex of the triangle, and the
Heisenberg picture on the left vertex of the triangle (see
also~\cite{Jacobs13a}). In these diagrams:
\begin{itemize}
\item The map $\Kl(\Rdn)\rightarrow \op{\EMod}$ on the left is the
  ``predicate'' functor, sending a space $X$ to the predicates on $X$,
  given by the effect module $\Cont(X,[0,1])$ of continuous functions
  $X\rightarrow [0,1]$, or for $C^*$-algebras mapping $A$ to the
  effects $[0,1]_{A}$. For $C^*$-algebras this was shown to be full
  and faithful in Lemma~\ref{CstarPredFunLem}, and for $\Kl(\Rdn)$ we
  combine Lemma~\ref{CstarPredFunLem} and Theorem~\ref{KlRToStarThm}:
$$\begin{array}{rcl}
\EMod\big(\Cont(Y,[0,1]), \Cont(X, [0,1])\big)
& = &
\EMod\big([0,1]_{C(Y)}, [0,1]_{C(X)}\big) \\
& \cong &
\HomPU\big(C(Y), C(X)\big) \\
& \cong &
\Kl(\Rdn)\big(X, Y\big).
\end{array}$$

\item The ``state'' functor $\Kl(\Rdn) \rightarrow \EM(\Rdn)$ is the
  standard full and faithful ``comparison'' functor from a Kleisli
  category to a category of Eilenberg-Moore algebras. In the
  $C^*$-algebra case it is the functor $\Stat$, combined with the
  equivalence from Theorem~\ref{SwirszczTheorem}. It is full and
  faithful by Theorem~\ref{StatFullFaithfulThm}.

\item The diagrams in~\eqref{KlRSETDiag} commute (up-to-isomorphism)
  in one direction. For $\Kl(\Rdn)$ we have:
$$\begin{array}{rcl}
\EMod\big(\Cont(X, [0,1]), [0,1]\big)
& = &
\EMod\big([0,1]_{C(X)}, [0,1]_{\C}\big) \\
& \cong &
\HomPU\big(C(X), \C)
\hspace*{\arraycolsep} = \hspace*{\arraycolsep}
\Rdn(X),
\end{array}$$
and similarly for $\CstarPU$ we have
\begin{align*}
\EMod([0,1]_A, [0,1]) & \cong \CstarPU(A, \C) 
   \qquad \mbox{by Lemma~\ref{CstarPredFunLem}} \\
 &= \Stat(A)
\end{align*}

\item The diagrams in~\eqref{KlRSETDiag} also commute (again,
  up-to-isomorphism) in the other direction, \textit{i.e.}
  $\Aff(\Rdn(X), [0,1]) \cong \Cont(X,[0,1])$ and
  $\Aff(\Stat(A),[0,1]) \cong [0,1]_A$. The former follows from the
  latter by taking $A = C(X)$, so we reduce to the latter. By Lemma
  \ref{ComplexKadisonLemma} we have that $A \cong \Aff(\Stat(A), \C)$
  as unital ordered vector spaces. We can then restrict both sides to
  their unit intervals and obtain an isomorphism $[0,1]_A \cong
  \Aff(\Stat(A), [0,1])$.
\end{itemize}\smallskip


\auxproof{
\begin{lem}
\label{IntervalIsAlgebraLemma}\hfill
\begin{enumerate}[(i)]
\item The unit interval $[0,1]$ is a $\Dst$-subalgebra of $(\C, \alpha)$. 
\item There is a $\Rdn$-algebra structure on $[0,1]$ such that the above $\Dst$-algebra structure is the underlying one constructed in lemma \ref{MonadFunctorLemma}.
\end{enumerate}
\end{lem}
\proof\hfill
\begin{enumerate}[(i)]
\item We need to show that, for $\varphi \in \Dst(\C)$, if $\supp \varphi \subseteq [0,1]$, $\alpha(\varphi) \in [0,1]$. This follows from the fact that if we have a finite set of numbers in $[0,1]$, their weighted mean lies in $[0,1]$. 

\auxproof
{
So assume we have such a $\varphi$. Then
\begin{eqnarray*}
\alpha(\varphi) & = & \sum_{z \in \C}\varphi(z) \cdot z \\
 & = & \sum_{z \in [0,1]}\varphi(z) \cdot z
\end{eqnarray*}
We prove this is in $[0,1]$ by proving the stronger statement that if $\min(\supp(\varphi)) \leq \alpha(\varphi) \leq \max(\supp(\varphi))$ by induction on $| \supp \varphi |$. 
\begin{itemize}
\item $| \supp \varphi | = 1$: Let $z$ be the unique $z \in \supp \varphi$, which is in $[0,1]$. We have tha $\varphi(z) = 1$ and hence $\alpha(\varphi) = z \in [0,1]$, and $z \leq z \leq z$.
\item Inductive step: Let $z'$ be the largest element of $\supp \varphi$. Then $\psi = \frac{\varphi|_{[0,1] \setminus \{z'\}}}{1 - \varphi(z')}$ is an element of $\Dst([0,1])$ to which we apply the inductive hypothesis and deduce that $\alpha(\psi) \in [0,1]$, and that $z' \geq \alpha(\psi)$. Since
\begin{eqnarray*}
\alpha(\varphi) & = & \sum_{z \in [0,1]} \varphi(z) \cdot z \\
 & = & \sum_{z \in [0,1] \setminus \{z'\}}\varphi(z) + \varphi(z') \cdot z' \\
 & = & (1-\varphi(z')) \sum_{z \in [0,1] \setminus \{z'\}}\frac{\varphi(z)}{1 - \varphi(z')} \cdot z + \varphi(z') \cdot z' \\
 & = & (1-\varphi(z'))\alpha(\psi) + \varphi(z')\cdot z'
\end{eqnarray*}
We show that the bounds apply to this as follows:
\begin{eqnarray*}
(1-\varphi(z'))\alpha(\psi) + \varphi(z')z' & \leq & (1-\varphi(z'))z' + \varphi(z')z' \\
 & = & z' 
\end{eqnarray*}
This establishes the upper bound as $z'$ was the largest element of $\supp \varphi$. For the lower bound we apply the induction hypothesis and get
\begin{eqnarray*}
\min(\supp(\varphi)) & \leq & \alpha(\psi) \\
 & = & (1 - \varphi(z'))\alpha(\psi) + \varphi(z')\alpha(\psi) \\
 & \leq & (1-\varphi(z'))\alpha(\psi) + \varphi(z')z' \\
 & = & \alpha(\varphi)
\end{eqnarray*}
\end{itemize}
}

\item Let $x \in C([0,1])$ be the function embedding $[0,1]$ into $\C$ in the usual way. We define an $\Rdn$-algebra structure on $[0,1]$ as follows:
\begin{eqnarray*}
\beta & \colon & \Rdn([0,1]) \rightarrow [0,1] \\
\beta(\phi) & = & \phi(x)
\end{eqnarray*}

We need the fact that $0 \leq \phi(x) \leq 1$ for all $\phi \in \Rdn([0,1])$. This true because any random variable taking values on the interval has its mean in the interval.
\auxproof{
Alternatively, proving this directly: Since $x$ is non-negative, by positivity of $\phi$ we have that $\phi(x) \geq 0$. To show the upper bound, we use that $1-x$ is non-negative, so $\phi(1-x) \geq 0$. Therefore:
\begin{eqnarray*}
\phi(1-x) = \phi(1) - \phi(x) & \geq & 0 \\
\phi(1) \geq \phi(x) \\
1 \geq \phi(x)
\end{eqnarray*}
}

So we now prove the unit law:
\begin{eqnarray*}
\beta(\eta_{[0,1]}(z)) & = & \eta_{[0,1]}(z)(x) \\
 & = & x(z) \\
 & = & z
\end{eqnarray*}

and the multiplication law:
\[
\xymatrix{
\Rdn^2([0,1]) \ar[r]^{\Rdn \beta} \ar[d]_{\mu_{[0,1]}} & \Rdn([0,1]) \ar[d]^\beta \\
\Rdn([0,1]) \ar[r]_{\beta} & [0,1]
}
\]
If $\phi \in \Rdn^2([0,1])$, then we have
\begin{itemize}
\item Bottom left: $ \beta(\mu_{[0,1]}(\phi)) = \mu_{[0,1]}(\phi)(x) = 
\phi(\lam{h}{h(x)})$.
\item Top right: $\beta(\Rdn(\beta)(\phi)) = \Rdn(\beta)(\phi)(x) = \phi(x \after \beta)$.
\end{itemize}
Then if $\psi \in \Rdn([0,1])$, we have that 
\begin{eqnarray*}
(x \after \beta)(\psi) & = & x(\beta(\psi)) \\
 & = & \beta(\psi) \\
 & = & \psi(x)
\end{eqnarray*}
Hence $x \after \beta = \lam{h}{h(x)}$, and so $\phi(x \after \beta) =
\phi(\lam{h}{h(x)})$, implying the commutativity of the diagram.

To see that $U\beta \after \tau_{[0,1]} = \alpha$, just observe that
integration commutes with finite summation (of measures).  \auxproof{
  Or directly:

Let $\varphi \in \Dst([0,1])$. Then
\begin{eqnarray*}
\beta(\tau_{[0,1]}(\varphi)) & = & \tau_{[0,1]}(\varphi)(x) \\
 & = & \sum_{z \in [0,1]} \varphi(z) \cdot x(z) \\
 & = & \sum_{z \in [0,1]} \varphi(z) \cdot z \\
 & = & \alpha(\varphi)
\end{eqnarray*}
}
\end{enumerate}
This establishes that if $([0,1],\beta)$ is an $\Rdn$-algebra, then $([0,1], \alpha)$ is the underlying $\Dst$-algebra. \qed
}

\noindent We summarise what we have just shown.

\begin{thm}
\label{SETThm}
The diagrams~\eqref{KlRSETDiag} are commuting ``state-and-effect'' triangles. \qed
\end{thm}

\section*{Final remarks}

The main contribution of this article lies in establishing a
connection between two different worlds, namely the world of
theoretical computer scientists using program language semantics (and
logic) via monads, and the world of mathematicians and theoretical
physicists using $C^*$-algebras. This connection involves the
distribution monad $\Dst$ on $\Sets$, which is heavily used for
modeling discrete probabilistic systems (Markov chains), in the
finite-dimensional case (see Proposition~\ref{FinKlDToStarFunProp})
and the less familiar Radon monad $\Rdn$ on compact Hausdorff spaces
(see Theorem~\ref{KlRToStarThm}). These results apply to both
commutative and noncommutative $C^*$-algebras, but only to positive
unital maps. Follow-up research will concentrate on characterizing
completely positive maps in the noncommutative case.



\subsubsection*{Acknowledgements} The authors wish to thank Hans Maassen, Jorik Mandemaker and Klaas Landsman for helpful discussions. 

This research has been financially supported by the Netherlands
Organisation for Scientific Research (NWO) under TOP-GO grant
no. 613.001.013 (The logic of composite quantum systems).

\bibliographystyle{plain}
\bibliography{gelfandjnl}

\begin{thebibliography}{10}

\bibitem{AbramskyC09}
S.~Abramsky and B.~Coecke.
\newblock A categorical semantics of quantum protocols.
\newblock In K.~Engesser, {Dov}~M. Gabbai, and D.~Lehmann, editors, {\em
  Handbook of Quantum Logic and Quantum Structures}, pages 261--323. North
  Holland, Elsevier, Computer Science Press, 2009.

\bibitem{Alfsen01}
E.M. Alfsen and F.W. Shultz.
\newblock {\em {State Spaces of Operator Algebras}}.
\newblock Birkh\"auser, 2001.

\bibitem{alfsen71}
Erik~M. Alfsen.
\newblock {\em {Compact Convex Sets and Boundary Integrals}}.
\newblock Ergebnisse der Mathematik und ihrer Grenzgebiete. Springer, 1971.

\bibitem{Alfsen80}
Erik~M. Alfsen, Harald Hanche-Olsen, and Frederic~W. Shultz.
\newblock {State Spaces of $C^*$-algebras}.
\newblock {\em Acta Mathematica}, 144(1):267--305, 1980.

\bibitem{Arveson81}
W.~Arveson.
\newblock {\em An Invitation to $C^*$-Algebra}.
\newblock Springer-Verlag, 1981.

\bibitem{Billingsley68}
Patrick Billingsley.
\newblock {\em {Convergence of Probability Measures}}.
\newblock John Wiley and Sons, 1968.

\bibitem{Conway90}
J.B. Conway.
\newblock {\em {A Course In Functional Analysis, Second Edition}}, volume~96 of
  {\em Graduate Texts in Mathematics}.
\newblock Springer Verlag, 1990.

\bibitem{Dixmier77}
J.~Dixmier.
\newblock {\em {$C^*$-Algebras}}, volume~15 of {\em North-Holland Mathematical
  Library}.
\newblock North-Holland Publishing Company, 1977.

\bibitem{DvurecenskijP00}
A.~Dvure\v{c}enskij and S.~Pulmannov{\'a}.
\newblock {\em New Trends in Quantum Structures}.
\newblock Kluwer Acad. Publ., Dordrecht, 2000.

\bibitem{FoulisB94}
D.~J. Foulis and M.K. Bennett.
\newblock Effect algebras and unsharp quantum logics.
\newblock {\em Found. Physics}, 24(10):1331--1352, 1994.

\bibitem{fremlin}
D.~H. Fremlin.
\newblock {Measure Theory, Volume 4}.
\newblock \url{http://www.essex.ac.uk/maths/people/fremlin/mt.htm}, 2003.

\bibitem{Fritz09}
Tobias Fritz.
\newblock {A Presentation of the Category of Stochastic Matrices}.
\newblock \url{http://arxiv.org/abs/0902.2554}, 2009.

\bibitem{Girard2011}
J-Y. Girard.
\newblock {Geometry of Interaction V: Logic in the hyperfinite factor}.
\newblock {\em Theor. Comput. Sci.}, 412(20):1860--1883, April 2011.

\bibitem{Giry82}
M.~Giry.
\newblock A categorical approach to probability theory.
\newblock In B.~Banaschewski, editor, {\em Categorical Aspects of Topology and
  Analysis}, volume 915 of {\em Lecture Notes in Mathematics}, pages 68--85.
  Springer Berlin Heidelberg, 1982.

\bibitem{Halmos50}
Paul~R. Halmos.
\newblock {\em {Measure Theory}}.
\newblock Number~18 in Graduate Texts in Mathematics. Springer, 1950.

\bibitem{HeinosaariZ12}
T.~Heinosaari and M.~Ziman.
\newblock {\em The Mathematical Language of Quantum Theory. From Uncertainty to
  Entanglement}.
\newblock Cambridge Univ. Press, 2012.

\bibitem{horodecki}
M.~Horodecki, P.~Horodecki, and R.~Horodecki.
\newblock {Separability of Mixed States: Necessary and Sufficient Conditions}.
\newblock {\em Physics Letters A}, 223(1–2):1 -- 8, 1996.

\bibitem{Jacobs12}
B.~Jacobs.
\newblock {\em Introduction to Coalgebra. Towards Mathematics of States and
  Observations}.
\newblock 2012.
\newblock Book, in preparation; version 2.0 available from
  \url{www.cs.ru.nl/B.Jacobs/CLG/JacobsCoalgebraIntro.pdf}.

\bibitem{Jacobs12f}
B.~Jacobs.
\newblock Involutive categories and monoids, with a {GNS}-correspondence.
\newblock {\em Found. of Physics}, 42(7):874--895, 2012.

\bibitem{Jacobs13a}
B.~Jacobs.
\newblock Measurable spaces and their effect logic.
\newblock In {\em Logic in Computer Science}. IEEE, Computer Science Press,
  2013.

\bibitem{JacobsM12b}
B.~Jacobs and J.~Mandemaker.
\newblock The expectation monad in quantum foundations.
\newblock In B.~Jacobs, P.~Selinger, and B.~Spitters, editors, {\em Quantum
  Physics and Logic (QPL) 2011}, volume~95 of {\em Elect. Proc. in Theor. Comp.
  Sci.}, pages 143--182, 2012.

\bibitem{JacobsM13a}
B.~Jacobs and J.~Mandemaker.
\newblock Relating operator spaces via adjunctions.
\newblock In J.~Chubb Reimann, V.~Harizanov, and A.~Eskandarian, editors, {\em
  Logic and Algebraic Structures in Quantum Computing and Information}, Lect.
  Notes in Logic. Cambridge Univ. Press, 2013.
\newblock See arxiv.org/abs/1201.1272.

\bibitem{Johnstone82}
P.~Johnstone.
\newblock {\em Stone Spaces}.
\newblock Number~3 in Cambridge Studies in Advanced Mathematics. Cambridge
  Univ. Press, 1982.

\bibitem{KadisonR83}
R.~Kadison and J.~Ringrose.
\newblock {\em Fundamentals of the Theory of Operator Algebras}.
\newblock Academic Press, 1983.

\bibitem{Kadison65}
Richard~V. Kadison.
\newblock {Transformations of States in Operator Theory and Dynamics}.
\newblock {\em Topology}, 3, Supplement 2(0):177 -- 198, 1965.

\bibitem{Kock70b}
A.~Kock.
\newblock On double dualization monads.
\newblock {\em Math. Scand.}, 27:151--165, 1970.

\bibitem{Maassen10}
H.~Maassen.
\newblock Quantum probability and quantum information theory.
\newblock In F.~Benatti, M.~Fannes, R.~Floreanini, and D.~Petritis, editors,
  {\em Quantum Information, Computation and Cryptography}, number 808 in Lect.
  Notes Physics, pages 65--108. Springer, Berlin, 2010.

\bibitem{Manes69}
E.~Manes.
\newblock A triple-theoretic construction of compact algebras.
\newblock In B.~Eckman, editor, {\em Seminar on Triples and Categorical Homolgy
  Theory}, number~80 in Lect. Notes Math., pages 91--118. Springer, Berlin,
  1969.

\bibitem{DualityMarkov}
M.~Mislove, J.~Ouaknine, D.~Pavlovic, and J.~Worrell.
\newblock {Duality for Labelled Markov Processes}.
\newblock In Igor Walukiewicz, editor, {\em Foundations of Software Science and
  Computation Structures}, volume 2987 of {\em Lecture Notes in Computer
  Science}, pages 393--407. Springer Berlin Heidelberg, 2004.

\bibitem{Mislove2012}
Michael Mislove.
\newblock {Probabilistic Monads, Domains and Classical Information}.
\newblock In Elham Kashefi, Jean Krivine, and Femke~van Raamsdonk, editors,
  {\em {\rm Proceedings 7th International Workshop on} Developments of
  Computational Methods, {\rm Zurich, Switzerland, 3rd July 2011}}, volume~88
  of {\em Electronic Proceedings in Theoretical Computer Science}, pages
  87--100. Open Publishing Association, 2012.

\bibitem{Moggi91a}
E.~Moggi.
\newblock Notions of computation and monads.
\newblock {\em Inf. \& Comp.}, 93(1):55--92, 1991.

\bibitem{Pelletier93}
J.~Wick Pelletier and J.~Rosick\'{y}.
\newblock {On the Equational Theory of $C^*$-algebras}.
\newblock {\em Algebra Universalis}, 30:275--284, 1993.

\bibitem{Rudin87}
W.~Rudin.
\newblock {\em Real and Complex Analysis}.
\newblock McGraw-Hill Book Company, 1987.
\newblock Third, International edition.

\bibitem{RussoDye66}
B.~Russo and H.A. Dye.
\newblock {A Note on Unitary Operators in $C\sp *$-algebras.}
\newblock {\em Duke Math. J.}, 33:413--416, 1966.

\bibitem{Sakai71}
S.~Sakai.
\newblock {\em $C^*$-algebras and $W^*$-algebras}, volume~60 of {\em Ergebnisse
  der Mathematik und ihrer Grenzgebiete}.
\newblock Springer, 1971.

\bibitem{Schaefer66}
Helmut~H. Schaefer.
\newblock {\em {Topological Vector Spaces}}, volume~3 of {\em Graduate Texts in
  Mathematics}.
\newblock Springer Verlag, 1966.

\bibitem{Semadeni73}
Z.~Semadeni.
\newblock {\em {Monads and their Eilenberg-Moore Algebras in Functional
  Analysis}}, volume~33 of {\em Queen's Papers in Pure and Applied
  Mathematics}.
\newblock Queen's University at Kingston, Ontario, Canada, 1973.

\bibitem{Street72}
R.~Street.
\newblock The formal theory of monads.
\newblock {\em Journ. of Pure \& Appl. Algebra}, 2:149--169, 1972.

\bibitem{Swirszcz74}
T.~\'Swirszcz.
\newblock {Monadic Functors and Convexity}.
\newblock {\em Bulletin de l'Acad\'{e}mie Polonaise des Sciences, S\'{e}rie des
  Sciences Math. Astr. et Phys.}, 22(1):39--42, 1974.

\bibitem{BramWesterbaan14}
A.~Westerbaan.
\newblock {Quantum Programs as Kleisli Maps}.
\newblock \texttt{http://arxiv.org/abs/1501.01020}.

\end{thebibliography}
\end{document}